\input amstex
\def\mapright#1{\smash{
   \mathop{\longrightarrow}\limits^{#1}}}
\def\mapleft#1{\smash{
   \mathop{\longleftarrow}\limits^{#1}}}
\def\mapdown#1{\Big\downarrow
   \rlap{$\vcenter{\hbox{$\scriptstyle#1$}}$}}

\def \cal{\Cal}
\documentstyle{amsppt}
\NoRunningHeads \magnification=\magstep1 \baselineskip=12pt
\parskip=5pt
\parindent=18pt
\topskip=10pt \leftskip=0pt \rightskip=0pt \pagewidth{30pc}
\pageheight{47pc}

\topmatter
\title Discriminant loci of ample and spanned line bundles
\endtitle
\author Antonio Lanteri and Roberto Mu\tildeaccent{n}oz \endauthor

\address Dipartimento di Matematica ``F. Enriques'',
Universit\`a, Via C. Saldini 50, I-20133 Milano, Italy
\endaddress
\email lanteri\@mat.unimi.it \endemail

\address
Departamento de Matem\'atica Aplicada, Universidad Rey Juan
Carlos, Calle Tulip\'an, 28933- M\'ostoles, Madrid, Spain
\endaddress
\email roberto.munoz\@urjc.es \endemail


\abstract Let $(X,L,V)$ be a triplet where $X$ is an irreducible
smooth complex projective variety, $L$ is an ample and spanned
line bundle on $X$ and $V\subseteq H^0(X,L)$ spans $L$. The
discriminant locus $\Cal D(X,V) \subset |V|$ is the algebraic
subset of singular elements of $|V|$. We study the components of
$\Cal D(X,V)$ in connection with the jumping sets of $(X,V)$,
generalizing the classical biduality theorem. We also deal with
the degree of the discriminant (codegree of $(X,L,V)$) giving some
bounds on it and classifying curves and surfaces of codegree 2 and
3. We exclude the possibility for the codegree to be 1.
Significant examples are provided.
\endabstract


\endtopmatter

\document

\head Introduction
\endhead

Let $X$ be an irreducible smooth complex projective variety of
dimension $n$. Take $L$ an ample line bundle on $X$ and a linear
system $|V| \subseteq |H^0(X,L)|$ with $\dim(|V|)=N$ and $V$
spanning $L$. We define the {\it discriminant locus} $\Cal D(X,V)$
of the triplet $(X,L,V)$ as the algebraic subset of $|V|$
parameterizing the singular elements of $|V|$. In the particular
case in which $\phi_V$ is an embedding, from now on the {\it
classical setting}, the discriminant locus is just the dual
variety $\phi_V(X)^\vee \subset \Bbb P^{N\vee}$, an irreducible
subvariety of $\Bbb P^{N\vee}$. A nice survey on results on
duality can be found in \cite{T}. When $\phi_V$ is not an
embedding some considerations on the morphism $\phi_V$ enter into
the picture. In fact the main ingredients to build $\Cal D(X,V)$
are the {\it jumping sets} (and their images by $\phi_V$),
measuring the deviation of $\phi_V$ from being an immersion, see
\cite{LPS1}. Inspired by the classical setting, different problems
on the discriminant locus can be faced. In our previous paper on
this subject \cite{LM1} (see also in \cite{LPS1}) we have focused
on the dimension of the discriminant locus. By the Bertini theorem
$\dim (\Cal D(X,V)) <N$. Hence it can be written as $\dim(\Cal
D(X,V))=N-1-k$, where $k \geq 0$ is called the {\it (discriminant)
defect} of $(X,L,V)$. Some bounds on $k$ and classification
results in the extremal cases (where $k$ is maximal) are provided
in \cite{LM1}. These results deeply rely on the geometry of
$\phi_V(X) \subset \Bbb P^N$ since $\phi_V(X)^\vee \subseteq \Cal
D(X,V)$. We have also studied this problem dropping the hypothesis
that $L$ is ample in \cite{LM2}.

When $\phi_V$ is an immersion $\phi_V(X)^\vee =\Cal D(X,V)$. The
locus where $\phi_V$ is not an immersion, consisting of the {\it
jumping sets}, is important to study the discriminant locus in
more general settings. In \cite{LPS1}, among other things, $\Cal
D(X,V)$ is written as a union of algebraic subsets built with the
{\it jumping sets} (see (0.3)). These sets are related with the
Chern classes of the first jet bundle of $L$. This approach is
continued in \cite{LPS2} where a partial study of this
decomposition of the discriminant locus (in the particular case in
which $\phi_V$ is generically one--to--one) is given. Some
considerations on the singular locus of a general $D \in \Cal
D(X,V)$ are also presented. In the current paper we follow this
line of research started in \cite{LPS1} and developed in
\cite{LPS2}, \cite{LM1}, \cite{LM2}. Our main goal is to find
appropriate generalizations of theorems holding in the classical
setting to the more general setting of an ample line bundle $L$
spanned by $V$.

A main theorem in the classical setting is the so called {\it
biduality theorem}. For $X \subset \Bbb P^N$ an irreducible
complex projective variety, $X^{\vee\vee}=X$ via the canonical
identification between $\Bbb P^N$ and $\Bbb P^{N\vee\vee}$. In
Section 1 we present a natural generalization of this theorem. In
fact we prove that any irreducible component of $\Cal D(X,V)$ is
the dual of the image of a component of a jumping set, see (1.3).
Moreover, the dual of any irreducible component of $\Cal D(X,V)$
is contained in $\phi_V(X)$ as proved in (1.4). These results help
to understand the relation between the decomposition in (0.3) and
the irreducible components of the discriminant. Significant
examples are provided.

Another basic fact in the classical setting is the irreducibility
of the dual variety of an irreducible complex projective variety.
In the non-classical setting this is no longer true. But if
$\phi_V$ is just an immersion, then $\Cal D(X,V)$ is still
irreducible. In Section 2 we show that, for curves, the facts of
$\phi_V$ being an immersion and the irreducibility of $\Cal
D(X,V)$ are equivalent. This is not true in higher dimension. We
can construct examples of surfaces for which the discriminant
locus is irreducible and any possible configuration of the
decomposition in (0.3) is achieved, $\phi_V$ not being, in
particular, an immersion. The most relevant consequence of
irreducibility of the discriminant locus is the emptiness of the
biggest jumping set, presented in (2.7).

Last problem we are concerned with is that of the degree of the
discriminant locus called, according to \cite{Z1}, {\it codegree
of $(X,L,V)$} (denoted codeg$(X,V)$). In the classical setting
this invariant is the {\it class} of $\phi_V(X) \subset \Bbb P^N$
when $\dim(\Cal D(X,V))=N-1$. In \cite{LPS1} it is shown that the
Chern classes of the first jet bundle are related with the
singular locus of elements in general linear subsystems of
appropriate dimension of $|V|$. Using this identification we get
an expression of the top Chern class of the first jet bundle
involving the degrees of the maximal dimensional components of the
discriminant. This expression and some consequences of it lead to
a complete classification of curves and surfaces of codegree less
than or equal to three. Let us recall that in the classical
setting a complete classification of smooth projective varieties
of codegree $\leq 3$ is provided in \cite{Z1}, \cite{Z2, Thm.
5.2}. We prove that there are no triplets $(X,L,V)$ with codegree
one and establish the complete list of curves and surfaces of
codegree two (see (4.4) and (6.11)) and three (see (4.4) and
(7.5)). All cases in the lists are effective and examples are
provided.

The final section is devoted to three further possible
developments of the theory. As a first thing we introduce the
concept of {\it tame codegree} for triplets $(X,L,V)$ for which
the general element in $\Cal D(X,V)$ is singular in just one point
and the singularity is quadratic and ordinary. This occurs in the
classical setting, but not only in this case. We classify (see
(8.1.4)) surfaces of tame codegree less than or equal to eight.
The second point is concerned with the study of the subvariety of
the discriminant made of the reducible or non-reduced elements in
$|V|$. The third one deals with two important facts holding in the
classical case for positive defect varieties but not yet explored
in the ample and spanned case: the parity theorem (the dimension
and the defect have the same parity) and the linearity of the
singular locus of a general element in the discriminant.

\head 0. Background material
\endhead

We work over the complex field and we use standard notation in
algebraic geometry. In particular, if $X$ is a projective
manifold, $K_X$ will denote the canonical bundle of $X$. We say
that a line bundle on $X$ is spanned by a vector space $V$ of
sections if $V$ generates $L$ at every point of $X$. By a little
abuse of notation line bundles and divisors are used with little
(or no) distinction. The symbol $\equiv$ denotes numerical
equivalence. We use the word {\it scroll} along the paper in the
classical sense, i.e., the projectivized of an ample vector bundle
with the polarization given by the tautological line bundle. We
fix our setting as follows.

\flushpar {\bf (0.0)} Let $(X,L,V)$ be a triplet where: $X$ is an
irreducible smooth projective variety of dimension $n$, $L$ is an
ample and spanned line bundle on $X$ and $V\subseteq H^0(X,L)$
spans $L$. Set $\dim(V)=N+1$ and let $\phi_V:X \to \Bbb P^N$ be
the morphism defined by $V$. In the particular case $V=H^0(X,L)$
we will write $\phi_L$.

The {\it{discriminant locus}} $\Cal D(X,V)$ of the triplet
$(X,L,V)$ parameterizes the singular elements of $|V|$. More
precisely, taking the incidence correspondence
$$\matrix
\Cal Y:=\{(x,[s]) \in X \times |V|:j_1(s)(x)=0\} & \mapright{p_1}
& X  \cr \mapdown{p_2} \cr \phantom{,,,,,,}\Cal D(X,V) \subset
\Bbb P^{N \vee},
\endmatrix \tag 0.1$$
where $j_1(s)$ denotes the first jet of the section $s \in V$,
$\Cal D(X,V)$ is the image of $\Cal Y$ via the second projection
of $X \times |V|$. Thus $\Cal D(X,V)$ is an algebraic subset in
$|V|=\Bbb P^{N \vee}$. By the Bertini Theorem $\dim (\Cal D(X,V))
< N$. Hence we can write $\dim (\Cal D(X,V)) = N-1-k$, where $k
\geq 0$ is called the {\it{defect}} of $(X,L,V)$. It is important
to point out the following fact.

\noindent {\bf (0.2)} We always look at the discriminant locus
$\Cal D(X,V) \subset |V|$ as an algebraic set with its reduced
structure.

If $\phi_V(X) \not= \Bbb P^N$ then the dual variety
$\phi_V(X)^{\vee}$ is a non-empty irreducible subvariety of $\Cal
D(X,V)$. Furthermore, if $\phi_V$ is an immersion then
$\phi_V(X)^{\vee} = \Cal D(X,V)$ \cite{LPS1, Rmk. 2.3.3}. Anyway,
points in $\Cal D(X,V) \setminus \phi_V(X)^{\vee}$ are coming from
points on $X$ where the differential of $\phi_V$ is not injective.
In this context it is natural to define the {\it jumping sets}
$\Cal J_i=\Cal J_i(V)$  ($1 \leq i \leq n$) as in \cite{LPS1,
(1.1)}, i.e., $\Cal J_i=\{x \in X: \text{rk}(d\phi_V(x)) \leq
n-i\}$. As in \cite{LPS2, (0.3.1)} $X_i$ stands for $\Cal J_i
\setminus \Cal J_{i+1}$, with the convention that $\Cal J_0=X$ and
$\Cal J_{n+1}=\emptyset$. This allows to define $\Cal D_i(X,V)
\subseteq \Cal D(X,V)$ as $\overline{p_2\circ p_1^{-1}(X_i)}$,
that is, the Zariski closure in $\Bbb P^{N \vee}$ of the locus of
elements of $|V|$ singular at points of $X_i$, so that: $$ \Cal
D(X,V) = \cup_{i=0}^n \Cal D_i(X,V). \tag 0.3$$

When no confusion arises we will write $\Cal D$ (respectively
$\Cal D_i$) instead of $\Cal D(X,V)$ (respectively $\Cal
D_i(X,V)$). For further use we  end this section with the
following easy consequence of the second Bertini theorem \cite{H,
III Ex. 11.3, p.\ 280}.

\proclaim{(0.4) Remark} For $(X,L,V)$ as in $(0.0)$, if $\dim(X)
\geq 2$ then any element of $|V|$ is connected. So if $D \in |V|$
is reducible then $D \in \Cal D(X,V)$.
\endproclaim

\head 1. On the components of the discriminant
\endhead

Let us study some properties of $\Cal D_i(X,V)$ ($0 \leq i \leq
n$) and relate them with the geometry of $\phi_V(X) \subseteq \Bbb
P^N$. A first basic fact is the following.

\flushpar {\bf (1.1)} $\Cal D_0(X,V)$ is always irreducible
because, if non-empty, it is the dual variety of $\phi_V(X)
\subset {\Bbb P}^N$.

This in fact does not mean that $\Cal D_0$ (when non-empty) is
always an irreducible component of $\Cal D$, as is shown, for
example, in \cite{LM1, Example 0.2}. Let us recall this example
for further references.

\noindent {\bf (1.1.0) Example.} Let $S$ be a Del Pezzo surface
with $K_S^2=1$ and let $L=-2K_S$. We know that $L$ is ample and
spanned and $\phi_L:S \to \Gamma \subset \Bbb P^3$ is a double
cover of the quadric cone $\Gamma$, branched at the vertex $v$ and
along the smooth curve $B$ cut out on $\Gamma$ by a transverse
cubic surface. We have $\Cal D(S,L) = \Cal D_0 \cup \Cal D_1 \cup
\Cal D_2$, where $\Cal D_0=\Gamma^\vee$ is a conic, the dual of
$\Gamma$, $\Cal D_1=B^\vee$ is the dual of $B$ and $\Cal
D_2=v^\vee$ is a plane. Recalling that $B$ is a sextic of genus
$4$ we thus get $\deg (\Cal D_1) = 2(\deg(B) + g(B) -1)=18$.
Furthermore we can note that $\Cal D_0 \subset \Cal D_1 \cap \cal
D_2$. Actually, $\Gamma^\vee \subset v^\vee$, since any plane
tangent to $\Gamma$ must contain $v$; moreover, $\Gamma^\vee
\subset B^\vee$ since any plane tangent to $\Gamma$ is tangent to
it along a generator $\ell$, hence it is also tangent to $B$ at
the points where $\ell$ meets $B$ (note that they are three
distinct points for the general $\ell$). On the other hand, note
that $B^\vee \cap v^\vee$ is a hyperplane section of $B^\vee$, so
it has degree $18$. Let $\ell_i$ be a generator of $\Gamma$
tangent to $B$. The line parameterizing the pencil of planes
through $\ell_i$ is contained in $B^\vee \cap v^\vee$. Actually
any plane in such pencil cuts $\Gamma$ along $\ell_i + \ell_i'$,
where $\ell_i'$ is another generator. So, this plane is in
$v^\vee$ (since containing $\ell_i$ it contains the vertex $v$);
moreover it is in $B^\vee$ (since $\ell_i$ is a line tangent to
$B$). Any such generator $\ell_i$ does correspond to a branch
point of the morphism $p: B \to \gamma$, where $\gamma =\Bbb P^1$
is a directrix of $\Gamma$. Since $p$ has degree $3$ and $B$ has
genus $4$, by Riemann--Hurwitz formula we get that this number is
$12$. All this shows that the intersection $\Cal D_1 \cap \Cal
D_2$ is given (scheme theoretically) by $3 \Cal D_0$ plus $12$
lines all tangent to $\Cal D_0$.

Assertion (1.1) is not true for $\Cal D_i(X,V)$ when $i>0$, as
shown by the following examples.

\flushpar {\bf (1.1.1) Examples.}

\noindent {\bf (a)} Consider the canonical system of a smooth
hyperelliptic curve of genus $g \geq 2$. The discriminant locus
consists of the union of $\Cal D_0$ and $\Cal D_1$, \cite{LPS1,
(1.8)}. In fact, $\Cal D_0=C^\vee$ is the dual variety of the
corresponding rational normal curve $C \subset \Bbb P^{g-1}$,
$\Cal D_1$ is reducible, being the union of $2g+2$ linear spaces
of dimension $g-2$, and $\Cal D_0 \setminus \Cal D_1\ne \emptyset
\ne \Cal D_1 \setminus \Cal D_0$.

\noindent {\bf (b)} Take $r>0$ triplets as in (0.0), say
$(X_1,L_1,V_1), \dots,(X_r,L_r,V_r)$, with the corresponding
morphisms $\phi_{V_i}:X_i \to \Bbb P^{N_i}$. Consider the product
morphism:$$ X=X_1 \times \dots \times X_r
\mapright{\phi_{V_1}\times \dots \times \phi_{V_r}} \Bbb P^{N_1}
\times \dots \times \Bbb P^{N_r}$$ and compose with the Segre
embedding to obtain $F:X \to \Bbb P^N$. For the triplet $(X,
L=F^*\cal O_{\Bbb P^N}(1),V=F^*H^0(\Bbb P^N,\Cal O_{\Bbb
P^N}(1)))$ it is straightforward to check the following fact:
$(x_1,\dots,x_r) \in \Cal J_i(X,V)$ if and only if $x_j \in \Cal
J_{i_j}(X_j,V_j)$ for $1 \leq j \leq r$ and $\Sigma_{j=1}^r (i_j)
\geq i$. Let us comment some particular cases.

(b.1) Take $r=2$, $X_1=C$, a smooth curve of genus $g$ and $|V_1|$
a base-point free pencil of degree $d$ defining a $d$--to--$1$ map
$\phi_{V_1}:C \to \Bbb P^1$. Choose $X_2=\Bbb P^{n-1}$, $L_2=\Cal
O_{\Bbb P^{n-1}}(1)$ and $V_2=H^0(X_2,L_2)$. For the triplet $(X,
F^*\cal O_{\Bbb P^N}(1),V=F^*H^0(\Bbb P^N,\Cal O_{\Bbb P^N}(1)))$,
$\Cal D(X,V)$ is the union of $\Cal D_0(X,V)= (\Bbb P^1 \times
\Bbb P^{n-1})^\vee$ (that is, $\Bbb P^1 \times \Bbb P^{n-1}
\subset \Bbb P^{2n-1 \vee} $) and $\Cal D_1(X,V)$, which is the
union of $s=2g-2+2d$ linear spaces of dimension
$N-1-(n-1)=2n-1-n=n-1$. Therefore, $\Cal D_0(X,V)=\Bbb P^1 \times
\Bbb P^{n-1} \subset \Bbb P^{2n-1} $, $\Cal D_1(X,V) =
\cup_{i=1}^s \Bbb P^{n-1}_i \subset \Cal D_0(X,V)$, and so $\Cal
D(X,V)=\Cal D_0(X,V)$.

(b.2) Now take $r=2$ and two triplets $(C_1,L_1,V_1)$ and
$(C_2,L_2,V_2)$ where $C_1$ and $C_2$ are smooth curves and, for
$i=1,2$, $|V_i|$ is a pencil of degree $d_i$. Consider the
corresponding morphisms $\phi_{V_i}:C_i \to \Bbb P^1$ and their
ramification loci $R_1=\{c_1,\dots,c_{s_{1}}\}\subset C_1$ and
$R_2=\{d_1,\dots,d_{s_{2}}\}\subset C_2$.  For the triplet $(X=C_1
\times C_2, F^*\Cal O_{\Bbb P^3}(1),V=F^*H^0(\Bbb P^3,\Cal O_{\Bbb
P^3}(1)))$ we have: $\Cal D_0(X,V)=(\Bbb P^1 \times \Bbb
P^1)^\vee=\Bbb P^1 \times \Bbb P^1 \subset \Bbb P^3$; $\Cal
J_1=\{(c,d) \in C_1\times C_2:c \in R_1\; \text{or}\; d \in R_2\}$
and $\Cal D_1(X,V) \subset \Cal D_0(X,V)$ is a union of lines;
$\Cal J_2=R_1 \times R_2$ and $\Cal D_2(X,V)$ is a union of
planes. Note that $\Cal D_2(X,V)$ is reducible and $\Cal
D(X,V)=\Cal D_0(X,V) \cup \Cal D_2(X,V)$.

(b.3) Let us recall here \cite{LPS2, Example 4.2.4}. Consider $C'
\subset \Bbb P^2$ an irreducible curve of degree $ \geq 4$ whose
singular locus is just a cusp. Call $\nu: C \to C'$ the
desingularization. Take $X_1=C$, $\phi_{V_1}$ the composition of
the desingularization with the inclusion $C' \subset \Bbb P^2$ and
$(X_2,L_2)=(\Bbb P^1, \Cal O_{\Bbb P^1}(1))$. In this situation
one can prove that $\Cal D=\Cal D_0$ since $\Cal D_1$ is a linear
space of dimension three contained in $\Cal D_0$.

\noindent {\bf (c)} Take a surface $\Sigma \subset \Bbb P^N$
having only an even set of nodes as singularities. One can take
the double cover $\pi: S \to \Sigma$, branched exactly at the
nodes. Here, ``even'' just means the following: consider the
blowing--up $Y \to \Sigma$ at the nodes, let $C_i$ be the
$(-2)$-curve corresponding to the node $p_i$ ($i=1, \dots ,\mu$),
and let $\Delta=\sum_{i=1}^{\mu}C_i$. The set of nodes of $\Sigma$
is even if $\Delta \in 2 \text{Pic}(Y)$. Under this condition, we
can consider the smooth surface $X$, double cover of $Y$ branched
along $\Delta$. Then the preimages on $X$ of the $C_i$'s are
$(-1)$-curves, and by contracting them we finally get the smooth
surface $S$ and the required double cover. Now let $L:=\pi^* \Cal
O_{\Sigma}(1)$ and $V= \pi^*W$, where $|W|$ is the trace of $|\Cal
O_{\Bbb P^N}(1)|$ on $\Sigma$. Then for our $(S,L)$, $\Cal J_2$
consists of $\mu$ points ($\mu$ being the number of nodes of
$\Sigma$). Moreover $\Cal J_1 \setminus \Cal J_2 = \emptyset$. So,
$\Cal D_2$ consists of $\mu$ hyperplanes, $\Cal D_1$ is empty and,
of course, $\Cal D_0$ is the dual of $\Sigma$.

This example is effective. Let $S=\text{J}C$ be the jacobian of a
smooth curve $C$ of genus $2$ and call $C$ again the image of the
curve embedded in $\text{J}C$ via the usual Abel--Jacobi map. Note
that $C$ is the theta divisor up to a translation, hence it is an
ample divisor. Set $L:=[2C]$. Then the ample line bundle $L$ is
also spanned, as Reider's theorem immediately shows; furthermore
$L^2=8$ and $h^0(L)=\chi(L)=L^2/2=4$. Moreover, $\phi_L: S \to
\Bbb P^3$ is a morphism of degree $2$ onto the Kummer quartic
surface $\Sigma$ having $16$ nodes as singular locus \cite{GH,
pp.\ 785--786}. This morphism of degree $2$ has exactly these $16$
points as branch locus, as can be checked by a local computation.
Then for this triplet $(S,L,H^0(S,L))$, $\Cal J_2$ consists of the
preimages of these $16$ points, while $\Cal J_1 \setminus \Cal
J_2= \emptyset$. Correspondingly, $\Cal D_2$ consists of $16$
planes in $\Bbb P^{3 \vee}=|V|$, and $\Cal D_1$ is empty. Note
also that $\Cal D_0 =(\Sigma)^\vee =\Sigma \subset \Bbb P^3$,
\cite{GH, p.\ 784}.

The following propositions generalize the fact that
$\phi_V(X)^\vee=\Cal D_0(X,V)$. Concretely, any irreducible
component of the discriminant locus is proved to be the dual
variety of the image by $\phi_V$ of an irreducible component of a
jumping set.

\proclaim{(1.2) Proposition} Let $(X,L,V)$ be a triplet as in
$(0.0)$ such that $\dim(X_i)=n-i$ and consider the irreducible
components of maximal dimension of $\overline{X}_i$, that is,
$Y_{ij} \subseteq \overline{X}_i$ ($1 \leq j \leq s_i$) such that
$\dim(Y_{ij})=n-i$. Then $\cup_{j=1}^{s_i}\phi_V(Y_{ij})^\vee
\subseteq \Cal D_i(X,V) \subseteq \Cal D(X,V).$
\endproclaim

\demo{Proof} Take a general point $y \in Y_{ij}$. Since $y \in
\Cal J_i \setminus \Cal J_{i+1}$ then $\text{rk}(d\phi_V(y))=n-i$.
Hence the kernel $K:=\text{ker}(d\phi_V(y))$ is a subspace of
dimension $i$ of the Zariski tangent space $T_{X,y}$. On the other
hand $\dim(Y_{ij})=n-i$. If $\dim(K \cap T_{Y_{ij},y})>0$ then
$\phi_V|_{Y_{ij}}$ is not finite, a contradiction. As a
consequence of the previous discussion one can choose local
coordinates $z_1, \dots, z_n$ around $y$ such that: (i) $\partial
s /\partial z_h=0$ for $1 \leq h \leq i$ and for all $s \in V$;
and (ii) $z_{i+1}, \dots, z_n$ are local parameters for $Y_{ij}$.
In this setting the vanishing of the derivatives with respect to
$z_{i+1}, \dots, z_n$ just means that the corresponding hyperplane
in $\Bbb P^N$ is tangent to $\phi_V(Y_{ij})$ at $\phi_V(y)$. \qed
\enddemo

Let us note that it can be $\dim(X_i)<n-i$ (for example in special
projections of smooth projective varieties). We can refer to
\cite{LPS2, Example 4.2.5} where a surface for which $\dim(X_1)=0$
and $X_2=\emptyset$ is provided. In fact it will be a consequence
of the next proposition that for \cite{LPS2, Example 4.2.5} $\Cal
D_1 \subset \Cal D_0$ and $\Cal D=\Cal D_0$.

\proclaim{(1.3) Proposition} Let $(X,L,V)$ be a triplet as in
$(0.0)$ and $\Cal D'$ an irreducible component of $\Cal D(X,V)$.
Then there exists an index $i$ ($0 \leq i \leq n$) and an
irreducible component $Y_{ij} \subseteq \overline{X}_i$ such that
$\dim(Y_{ij})=n-i$ and $\Cal D'=\phi_V(Y_{ij})^\vee.$
\endproclaim

\demo{Proof} Consider the following incidence correspondence
$$\matrix
\Cal Y_{|\Cal D'}:=\{(x,[s]) \in X \times \Cal D':j_1(s)(x)=0\} &
\mapright{p_1} & X  \cr \mapdown{p_2} \cr \phantom{,,,,,,}\Cal D'
\subset \Bbb P^{N \vee}.
\endmatrix$$
Let $\dim (\Cal D')=N-1-k'$. By \cite{LPS2, Lemma (0.6)} the
dimension of the generic fibre of $p_2$ is $k'=N-1-\dim(\Cal D')$
and so $\dim(\Cal Y_{|\Cal D'})=N-1$. Take a $(N-1)$--dimensional
irreducible component $\Cal Y^0 \subseteq \Cal Y_{|\Cal D'}$ such
that $p_2(\Cal Y^0)=\Cal D'$. Let $i$ be the maximum integer such
that $p_1(\Cal Y^0) \subseteq \overline{X}_i$, and consider a
general $p \in p_1(\Cal Y^0) \cap X_i$. As $p \in X_i$ it holds
that $|V-2p| \subseteq \Cal D(X,V)$ is a linear space $T_p$ of
dimension $N-1-(n-i)$. Since $T_p \cap \Cal D' \ne \emptyset$ then
$T_p \subseteq \Cal D'$. Whence the dimension of the general fibre
of $p_1$ is $N-1-(n-i)$. In particular, $\dim(X_i)=n-i$ and there
exists an irreducible component $Y_{ij} \subseteq \overline{X}_i$
such that $Y_{ij}=p_1(\Cal Y^0)$. Just by definition of dual
variety, $\phi_V(Y_{ij})^\vee =\Cal D'$. \qed
\enddemo

\noindent (1.3.1) With the same notation as in the proof of (1.3)
we have maps $\Cal D' \mapleft {p_2} \Cal Y^0 \mapright{p_1} X$
where $\Cal Y^0$ is characterized by the following properties:
irreducibility, $\dim(\Cal Y^0)=N-1$ and $p_2(\Cal Y^0)=\Cal D'$.
Moreover, by classical biduality theorem, any other
$(N-1)$--dimensional irreducible component $\Cal Y^1 \subseteq
\Cal Y_{|\Cal D'}$ such that $p_2(\Cal Y^1)=\Cal D'$ verifies
$\phi_V(p_1(\Cal Y^1))=\phi_V(Y_{ij})$.

As a consequence we have the following statement, analogous to the
classical {\it biduality theorem}, offering in particular some
control on the linear components of $\Cal D(X,V)$.

\proclaim{(1.4) Biduality Theorem} Let $(X,L,V)$ be a triplet as
in $(0.0)$. Then $(\Cal D')^\vee \subseteq \phi_V(X)$ for any
irreducible component $\Cal D' \subseteq \Cal D(X,V)$.
\endproclaim

\demo{Proof} From (1.3) and with the same notation as there it
holds that $\Cal D'=\phi_V(Y_{ij})^\vee$. Then our result is a
consequence of the classical biduality theorem, that is, $(\Cal
D')^\vee =(\phi_V(Y_{ij})^\vee)^\vee=\phi_V(Y_{ij})\subseteq
\phi_V(\overline{X_i})\subseteq \phi_V(X)$. \qed
\enddemo

\flushpar (1.4.1) Let $(X,L,V)$ be a triplet as in $(0.0)$.
Suppose there exists a linear irreducible component $\Bbb
P^{N-1-r} \subseteq \Cal D(X,V)$. By biduality $\phi_V(X)$
contains a linear space $T$ of dimension $r$ and by (1.3) there
exists $Y \subset \Cal J_{n-r}$ such that $T=\phi_V(Y)$. It is of
particular interest the fact that if $\Cal D(X,V)$ contains a
hyperplane, then $\Cal J_n \ne \emptyset$. In fact, any hyperplane
in $\Cal D(X,V)$ defines a point in $\Cal J_n$. Note that the
converse is obvious, because $|V-x| \subset \Cal D(X,V)$ if $x \in
\Cal J_n$. So we have the following

\proclaim{(1.4.2) Corollary} Let $(X,L,V)$ be a triplet as in
$(0.0)$. Then $\Cal D(X,V)$ contains a hyperplane if and only if
$\Cal J_n \not= \emptyset$.
\endproclaim

\proclaim{(1.5) Lemma} Let $(X,L,V)$ be a triplet as in $(0.0)$
and let $\dim(\Cal D(X,V))=N-1-k$. Then $\Cal J_i =\emptyset$ for
$i>n-k$. Moreover, if $\dim(X_{n-k})=k$ then any maximal
dimensional component of $\overline{\phi_V(X_{n-k})}$ is linear.
\endproclaim

\demo{Proof} Let us suppose there exists $i >n-k$ for which $X_i
\ne \emptyset$. For any $p \in X_i$ we get $|V-2p|=\Bbb
P^{N-1-(n-i)} \subseteq \Cal D(X,V)$ a contradiction. The last
assertion is just a consequence of (1.3). \qed \enddemo

Let $(X,L,V)$ and $(Y,M,W)$ be two triplets as in (0.0) such that
$\dim(V)=\dim (W)=N+1$. In the classical case, that is, $\phi_{V}$
and $\phi_{W}$ embeddings, the {\it biduality theorem} states that
if $\Cal D(X,V)=\Cal D(Y,W)\subset \Bbb P^N$ (that is, there
exists a linear transformation of $\Bbb P^N$ sending
isomorphically $\Cal D(X,V)$ to $\Cal D(Y,W)$) then $(X,L)=(Y,M)$.
It is natural to ask to what extent this theorem is true when
$\phi_V$ or $\phi_W$ are not embeddings. Next examples show that
it cannot be true in the same terms and the right hypotheses to
impose.

\noindent {\bf (1.6) Examples.}

{\bf (a)} Choose $X_1$ a smooth elliptic curve and $L_1$ giving a
$g^1_2=|V_1|$ on $X_1$. Take $(X_2,L_2,V_2)=(\Bbb P^1, \Cal
O_{\Bbb P^1}(1),H^0(\Bbb P^1, \Cal O_{\Bbb P^1}(1)))$. As in (b)
of (1.1.1) we have $(X,L,V)$ such that $\phi_V(X) \subset \Bbb
P^3$ is a smooth quadric and the branch locus of $\phi_V$ is made
of four disjoint lines on $\phi_V(X)$. Hence $\Cal D(X,V)$ is a
smooth quadric in $\Bbb P^3$. For $(Y=Q,L=\Cal O_Q(1),V=H^0(Y,L))$
a smooth quadric with its corresponding embedding in $\Bbb P^3$ we
have $\Cal D(X,V)=\Cal D(Y,W)$, $\phi_V(X)=\phi_W(Y) \subset \Bbb
P^3$ but $X$ and $Y$ are not isomorphic.

{\bf (b)} Choose $X_1$ a smooth conic such that
$(X_1,L_1,V_1=H^0(X_1,L_1))$ defines the embedding
$\phi_{V_1}(X_1) \subset \Bbb P^2$. Take $X_2$ a double cover of
the plane, $f:X_2 \to \Bbb P^2$, branched along $\phi_{V_1}(X_1)$,
$L_2=f^*\Cal O_{\Bbb P^2}(1)$, $V_2=f^*H^0(\Bbb P^2,\Cal O_{\Bbb
P^2}(1))$. Then $\phi_{V_1}(X_1)\subset \Bbb P^2$ is a smooth
conic and $\phi_{V_2}(X_2)=\Bbb P^2$. In fact $\phi_{V_1}(X_1)$
and $\phi_{V_2}(X_2)$ are not isomorphic but $\Cal
D(X_1,H^0(X_1,L_1))=\Cal D(X_2,V_2)=\phi_{V_1}(X_1)^\vee$ a smooth
conic.

{\bf (c)} Consider two smooth plane curves $C_1,C_2 \subset \Bbb
P^2$, not isomorphic. Let $f_1: X_1 \to C_1 \times \Bbb P^2$ and
$f_2: X_2 \to \Bbb P^2 \times C_2$ be cyclic covers, both branched
along $C_1 \times C_2$. Let $F_i$ be the composition of the Segre
embeding $\Bbb P^2 \times \Bbb P^2 \subset \Bbb P^8$ with $f_i$.
Then we get triplets $(X_i, L_i=F_i^*\Cal O_{\Bbb P^8}(1),
V_i=F_i^*H^0(\Bbb P^8,\Cal O_{\Bbb P^8}(1)))$, $i=1,2$. Whence
$\Cal D(X_1,V_1)= \Cal D_0 \cup \Cal D_1$ where $\Cal D_0 =(C_1
\times \Bbb P^2)^\vee$ and $\Cal D_1=(C_1 \times C_2)^\vee$. We
know $\dim(\Cal D_0)=6$ (since $\Cal D_0$ is the dual of a
three-dimensional scroll over a curve). We claim that $\Cal D_0
\subset \Cal D_1$. In fact a general element of $\Cal D_0$
corresponds to a hyperplane $H \in (C_1 \times \Bbb P^2)^\vee$
which is tangent to $C_1 \times \Bbb P^2$ along a line contained
in a fiber $f$. Since this line is meeting $f \cap (\Bbb P^2
\times C_2)=C_2$ then $H \in (C_1 \times C_2)^\vee$. Hence $\Cal
D(X_1,V_1)=(C_1 \times C_2)^\vee$. In the same way we see that
$\Cal D(X_2,V_2)=(C_1 \times C_2)^\vee$. Note however that
$\phi_{V_1}(X_1)=C_1 \times \Bbb P^2$ is not isomorphic to
$\phi_{V_2}(X_2)=\Bbb P^2 \times C_2$.

\proclaim{(1.7) Proposition} Let $(X,L,V)$ and $(Y,M,W)$ be two
triplets as in $(0.0)$ such that $\dim(X)=n=\dim(Y)$ and
$\dim(V)=\dim (W)=N+1$. If $\Cal D_0(X,V)$ is an irreducible
component of $\Cal D(Y,W)$ then $\phi_V(X)=\phi_W(Y) \subset \Bbb
P^N$.
\endproclaim

\demo{Proof} Since $\Cal D_0(X,V)=\phi_V(X)^\vee$ is a component
of $\Cal D(Y,W)$ then, by (1.3), there exists $Y_{ij} \subset Y$
such that $\phi_V(X)^\vee =\phi_V(Y_{ij})^\vee$. By the classical
biduality theorem $\phi_V(X) =\phi_V(Y_{ij})$. This gives
$\dim(\phi_V(Y_{ij}))=\dim (Y)$ and so $Y_{ij}=Y$. \qed
\enddemo

\head 2. Irreducibility of the discriminant locus
\endhead

In this section we study some consequences of the irreducibility
of the discriminant locus. A general fact is the following.

\flushpar {\bf (2.1)} If $\phi_V$ is an immersion then $\Cal
D(X,V)=\Cal D_0$ and so it is irreducible, see (1.1).

The converse of (2.1) is also true for curves. We need the
following lemma.

\proclaim{(2.2) Lemma} Let $(X,L,V)$ be as in $(0.0)$ such that
$\Cal D(X,V)$ is irreducible and $\Cal J_n \ne \emptyset$, then:

$(2.2.1)$ $\Cal J_n$ is a finite set and $\phi_V(x)=\phi_V(y)$ for
any $x,y \in \Cal J_n$;

$(2.2.2)$ $\Cal D(X,V)$ is a hyperplane of $|V|$ and $\phi_V(X)$
is a cone whose vertex contains $\phi_V(\Cal J_n)$.
\endproclaim

\demo{Proof} By \cite{LPS1, Theorem 1.2} $\dim(\Cal J_n)=0$. For
any $x \in \Cal J_n$ one has $\Bbb P^{N-1}=|V-x| \subseteq \Cal
D(X,V)$. Since $\Cal D(X,V)$ is irreducible then $|V-x|=|V-y|$ for
any $x,y \in \Cal J_n$ and (2.2.1) follows. Moreover $\Cal
D(X,V)=\Cal D_n=|V-x|$ for any $x \in \Cal J_n$. Since
$\phi_V(X)^{\vee}=\Cal D_0 \subseteq \Cal D(X,V)$ then it is
either empty or non-empty and degenerate (in the sense that it is
contained in a hyperplane of $\Bbb P^N$). If empty then
$\phi_V(X)=\Bbb P^N$ and $N=n$. So, $\phi_V(X)=\Bbb P^n$ is a
linear cone. If non-empty and degenerate then $\phi_V(X) \subset
\Bbb P^N$ is a cone whose vertex contains $\phi_V(x)$ for any $x
\in \Cal J_n$. \qed
\enddemo

We will see in (3.5) that (2.2.2) cannot occur.

\flushpar{{\bf (2.3) Remark.} For $C$ a smooth irreducible curve,
it is not possible to construct a finite morphism $\pi:C \to \Bbb
P^1$ of degree $d \geq 2$ with a single branch point $p \in \Bbb
P^1$. Let us call $m$ the number of distinct points in
$\pi^{-1}(p)$. The claim is just a consequence of the
Riemann--Hurwitz formula: $2g(C)-2=-2d+d-m$.

We can prove the following result for curves.

\proclaim{(2.4) Proposition} Let $(C,L,V)$ be as in $(0.0)$ with
$\dim(C)=1$. Then $\Cal D(C,V)$ is irreducible if and only if
$\phi_V$ is an immersion.\endproclaim

\demo{Proof} In view of (2.1) it is enough to prove the {\it only
if} part. So let us assume $\Cal D(C,V)$ to be irreducible. If
$\Cal J_1 \ne \emptyset$ then $\phi_V(C)=\Bbb P^1$ by (2.2.2). The
contradiction comes from (2.2.1) because $\Cal J_1 \ne \emptyset$
implies that one can construct a map as in (2.3). \qed
\enddemo

This statement is not true for higher dimension. In the examples
(b.1) and (b.3) of (1.1.1) (see \cite{LPS2, Example 4.2.4}) $\Cal
D(X,V)$ is not only irreducible but equal to $\phi_V(X)^\vee$ not
being $\phi_V$ an immersion. Let us record a list of examples of
surfaces with irreducible discriminant locus showing different
behaviors of the $\Cal D_i's$. Since scrolls are of particular
interest, first consider the following result, that is important
also for the next sections. We follow the usual notation of
\cite{Ha, V 2}.

\proclaim{(2.5) Lemma} Let $(S,L,V)$ be a triplet as in $(0.0)$
such that $\dim(S)=2$ and $(S,L)$ is a scroll over a smooth curve
$B$. Let $C_0$ and $f$ be a fundamental section and a fibre
respectively. Then

$(2.5.1)$ $L \equiv C_0+bf$, with $b > 0$.

$(2.5.2)$ $\Cal J_2=\emptyset$.

$(2.5.3)$ If $\dim(V)=3$ there is a surjection $i:B \to \Cal
D(S,V)$. In particular $i$ is an isomorphism if $b=1$.
\endproclaim

\demo{Proof} Write $S = \Bbb P(\Cal E)$, where $\Cal E$ is a
vector bundle of rank $2$ on $B$. We can assume that $\Cal E$ is
normalized in the sense of \cite{Ha, p.\ 373} and that $C_0$ is
the tautological section on $S$. Of course $L \equiv C_0+bf$ for
some integer $b$, since $(S,L)$ is a scroll. Let $\pi:S \to B$ be
the scroll projection. Since $L$ is ample and the general element
in $|L|$ is irreducible then $b \geq 0 $ by \cite{Ha, V Prop. 2.20
and 2.21}. It is not hard to see that $b=0$ implies that $L$ is
not globally generated and this proves (2.5.1). Alternatively an
argument based on the non emptiness of the discriminant locus can
be given. Take $D \in \Cal D(S,V)$ ($\Cal D(S,V) \not= \emptyset$
by \cite{LPS, Thm. 2.8}). Choose $x \in \text{Sing}(D)$, and let
$f_{\pi(x)}=\pi^*\Cal O_B(\pi(x))$ be the fibre of $S$ containing
$x$. Then
$$D = f_{\pi(x)} + R \tag 2.5.4$$ for some effective divisor $R
\equiv C_0 +(b-1)f$ containing $x$. Otherwise $1 = L f = D
 f_{\pi(x)} \geq \text{mult}_x(D)
\text{mult}_x(f_{\pi(x)}) \geq 2,$ a contradiction. Suppose that
$b \leq 0$. This would mean that $0 < h^0(R) = h^0(\Cal E \otimes
\Cal L)$ where $\Cal L$ is a line bundle on $B$ with $\deg( \Cal
L) <0$, contradicting the assumption that $\Cal E$ is normalized
made at the beginning.

To prove (2.5.2), assume that $x \in \Cal J_2$. Then
$|V-x|=|V-2x|$. As we have seen before any $D\in |V-2x|$ is as in
(2.5.4). Hence $|V-2x|=f_{\pi(x)}+|V-f_{\pi(x)}-x|$. Then
$f_{\pi(x)}$ would be contained in the base locus of $|V-x|$. But
$\text{Bs}(|V-x|)=\phi_V^{-1}(\phi_V(x))$ must be a finite set
since $L$ is ample and spanned by $V$. This gives a contradiction.

Now assume that $\dim(V) = 3$, so that $\phi_V(S)=\Bbb P^2$. For
every $p \in B$ let $x, y$ be any two distinct points lying on the
fibre $f_p=\pi^{-1}(p)$. So $|V-x-y|$ consists of a single element
$D_p$, and $D_p=f_p+R_p$ for an effective $R_p$ such that $R_p
\equiv C_0+(b-1)f$. In particular, $D_p$ has a (exactly one)
singular point on $f_p$. Then the mapping $i(p)= D_p$ defines a
morphism $i:B \to \Cal D(S,V)$. Now pick an element $D \in \Cal
D(S,V)$ and let $x$ be a singular point of $D$. By (2.5.4) it
holds that $D=i(\pi(x))$. Then $i$ is surjective. If $b=1$ then
$i$ is also injective. If $i(p)=i(q)$ with $p \not=q$, then
$D_p=D_q$; in particular $D_q-f_p-f_q \equiv C_0-f$ would be
effective, a contradiction.
 \qed
\enddemo

In the previous lemma we have shown that for any $N=\dim(|V|)$ we
have the following maps:
$$\matrix \Cal B:=\{(b,D) \in B \times \Cal D(S,V):\pi^{-1}(b)
\subset D \} & \mapright{\pi_1} & B \cr \mapdown{\pi_2} \cr
\phantom{,,,,,,}\Cal D(S,V) \subset \Bbb P^{N \vee},
\endmatrix $$ where any fibre of $\pi_1$ is a linear
space of dimension $N-2$. If (2.5.3) holds, $\pi_1$ is an
isomorphism since $N=2$ and $i=\pi_2\circ \pi_1^{-1}$.

We are going to show several examples of surfaces whose
discriminant locus is irreducible. The list shows that any
possible relation between $\Cal D_0$, $\Cal D_1$ and $\Cal D$ can
occur. In fact (2.4) is no longer true when the dimension is
bigger than one. In the following examples $\Cal D_2=\emptyset$,
being $\Cal D$ irreducible. This is a general fact, as it will be
proved later, see (2.7).

\flushpar{\bf (2.6) Examples.}

{\bf (a)} Take $C \subset \Bbb P^3$ an irreducible non-degenerate
smooth curve. As in \cite{LM, Example 3.2} consider the conormal
variety $\Cal X=\{(c,H): T_{C,c} \subset H\} \subset C \times
C^\vee$ and the corresponding projections $\pi_1$ and $\pi_2$. The
triplet $$(\Cal X, \pi_2^*\Cal O_{\Bbb P^{3\vee}}(1),
\pi_2^*H^0(\Bbb P^{3 \vee},\Cal O_{\Bbb P^{3\vee}}(1)))$$ is as in
$(0.0)$. Recall that $\Cal X$ is a $\Bbb P^1$-bundle over $C$. A
local computation shows that $\Cal J_1$ consists of a section plus
the fibres over the hyperflexes of $C$. In fact the section
corresponds to $\overline{\{(c,\text{Osc}_c^2(C)):\; c \in U}\}
\subset C \times C^{\vee}$, where $\text{Osc}_c^2(C)$ stands for
the second osculating projective linear space to $C$ at $c \in C$
and $U \subset C$ is the open subset of points of $C$ where the
osculating plane is defined. One can check that $\Cal D$ consists
of the tangent developable $TC$ of $C$ plus the lines
corresponding to the hyperflexes. Then $\Cal D_0=(C^\vee)^\vee=C
\subset \Cal D = \Cal D_1$.

{\bf (b)} Cyclic coverings of $\Bbb P^2$ give rise to:
$\emptyset=\Cal D_0 \subset \Cal D_1=\Cal D=B^\vee$ where $B$ is
the branch locus of the covering.

{\bf (c)} In the example (b.3) of (1.1.1), see \cite{LPS2, Example
4.2.4}, we have $\emptyset \ne \Cal D_1 \subset \Cal D_0=\Cal D$.
If $\phi_V$ is an immersion then $\emptyset = \Cal D_1 \subset
\Cal D_0=\Cal D$.

{\bf (d)} Let $B$ be a smooth elliptic curve and $p \in B$.
Consider the rank two vector bundle $\cal E$ defined as the
non-trivial extension: $0 \to \Cal O_B \to \Cal E \to \cal O_B(p)
\to 0.$ Set $S=\Bbb P(\Cal E)$, $L\equiv C_0+f$ and $V=H^0(S,L)$.
One can check that $h^0(S,L)=3$ and $L$ is ample and spanned. By
(2.5) there is an isomorphism between $B$ and $\Cal D(S,V)$. Then
$\phi_V:S \to \Bbb P^2$ is a degree 3 map whose branch locus is
the dual of the smooth plane cubic $\Cal D(X,V)$. In this case
$\emptyset =\Cal D_0 \subset \Cal D_1 =\Cal D$.

The remainder of this section is devoted to prove that
irreducibility of the discriminant locus implies emptiness of the
maximal jumping set. As (2.6) shows in the case of surfaces, this
is in principle the only consequence of irreducibility of the
discriminant locus one can expect in general.

\proclaim{(2.7) Lemma} Let $(X,L,V)$ be as in $(0.0)$. If $\Cal
D(X,V)$ is irreducible then $\Cal J_n = \emptyset$.
\endproclaim

\demo{Proof} If $\Cal J_n \ne \emptyset$ and $\Cal D(X,V)$ is
irreducible then $\Cal D(X,V)=\Bbb P^{N-1}=|V-x|$ for any $x \in
\Cal J_n$. Consider a general $W \subseteq V$ such that
$\dim(W)=n+1$ and note that $W$ spans $L$. Then $\Cal
D(X,W)=|W-x|$ and $\phi_W(X)=\Bbb P^n$. In particular there exists
$p \in \phi_W(X)$ such that $p=\phi_W(x)$ for any $x \in \Cal
J_n$. Moreover, for this choice of $W$, $\Cal J_1=\Cal J_1(W)$ is
a divisor on $X$ and for any component $\Cal J_1^{i} \subseteq
\Cal J_1$ we have $\phi_W(\Cal J_1^{i})^\vee \subset \Cal D(X,W)$.
That is, $\phi_W(\Cal J_1^{i})^\vee$ is contained in the
hyperplane of $\Bbb P^{N \vee}$ of hyperplanes of $\Bbb P^N$
through $p$. Then

\noindent (2.7.1) $\phi_W(\Cal J_1)$ is a union of cones with
vertex containing $p$.

\noindent Choose a general line $T \subset \Bbb P^n$ through $p$.
Since $T \cap \phi_W(\Cal J_1)=\{p\}$ then, by Bertini type
theorems, $\phi_W^{-1}(T)$ is a curve whose singular locus is
contained in $\phi_W^{-1}(p)$. In fact, consider a general
hyperplane $H \subset \Bbb P^n$ and define $f:X \setminus \Cal J_1
\to H \setminus \phi_V(\Cal J_1)$ as $f(z)=\langle \phi_W(z),
\phi_W(p) \rangle \cap H$. Hence $T \setminus \phi_W^{-1}(p)$ is
smooth because it is a general fibre of $f$. Then, consider any
component $\Gamma_i \subset \phi_W^{-1}(T)$ and let
$\mu_i:\gamma_i \to \Gamma_i$ be its desingularization. The
morphism $\phi_W \circ \mu: \gamma_i \to \Bbb P^1$ has only one
branch point, so, by (2.3), it is an isomorphism. This says that
$\Gamma_i$ is a smooth rational curve such that $L \Gamma_i=1$.
This implies in particular that $X$ is swept out by lines, that
is, there exists a family $\Cal T$ of smooth rational curves of
$L$--degree $1$, sweeping out $X$. Moreover, since
$\phi_V^{-1}(p)$ is finite, there exists $x \in \phi_V^{-1}(p)$
such that $x \in \ell$ for $\ell$ general in $\Cal T$. By
\cite{LPS2, Lemma 3.1} the normal bundle $N_{\ell/X}$ splits as
$\Cal O_{\Bbb P^1}(a_1)\oplus \dots \oplus \Cal O_{\Bbb
P^1}(a_{n-1})$ with $0 \leq a_1\leq \dots \leq a_{n-1}$. Let us
suppose that $a_1= \dots =a_t=0,$ $a_{t+1}>0$. Since
$\phi_W(X)=\Bbb P^n$ and $x \in \ell$ for $\ell \in \Cal T$
general we conclude that the irreducible component $\Cal T_x$
through $\ell$ of the Hilbert scheme of rational curves of
$L$--degree 1 on $X$ through $x$ has dimension greater than or
equal to $n-1$. Since $h^1(\ell,N_{\ell/X}(-1))=0$ then $\Cal T_x$
is smooth; moreover $h^0(\ell,N_{\ell/X}(-1))=a_{t+1}+ \dots
+a_{n-1} \geq n-1$. Then $h^0(\ell,N_{\ell/X}) \geq n-1+n-1=2n-2$.
This implies $(X,L)=(\Bbb P^n,\Cal O_{\Bbb P^n}(1))$ by \cite{LP2,
Thm. 1.4}, a contradiction. \qed
\enddemo

\head 3. Codegree  \endhead

The definition of {\it codegree} can be established as in the
classical case, see \cite{Z1}.

\proclaim{(3.1) Definition} For $(X,L,V)$ as in $(0.0)$ define its
{\rm codegree}, say {\rm codeg}$(X,V)$, as the degree of $\Cal
D(X,V) \subset \Bbb P^{N\vee}$.
\endproclaim

As said in (0.2), $\Cal D(X,V)$ is an algebraic subset of $|V|$
with its reduced structure; hence codeg$(X,V)=\Sigma_{\dim(\Cal
D^{i})={\text max}}$deg$(\Cal D^{i})$, the sum ranging over all
irreducible components of maximal dimension.

As in the classical case we can relate the codegree of $(X,L,V)$
with the Chern classes of the jet bundle $J_1(L)$. Suppose
$\dim(\Cal D(X,V))=N-1$ and consider the maximal dimensional
components of its discriminant locus, say $\cal D^1, \dots, \Cal
D^s \subseteq \Cal D(X,V)$. Let $d_i=\text{deg}(\Cal D^{i})$. For
$D \in \Cal D^{i}$, $1 \leq i \leq s$, with isolated singularities
take $s \in V$ the section defining $D$ and $x \in $Sing$(D)$.
Define the 0--cycle $z_x=\mu_x(D)x$, where $\mu_x(D)$ is the
Milnor number of $x$ as an isolated singular point of $D$, and
$z(D)=\Sigma_{x \in {\text Sing}(D)}z_x$. Then we can prove the
following

\proclaim{(3.2) Theorem} Let $(X,L,V)$ be as in $(0.0)$ with $\dim
(\Cal D)=N-1$ and consider a general pencil $|W|\subseteq |V|$.
With the notation of the previous paragraph, let $|W| \cap \Cal
D^{i} = \{D^{i}_j\},\ (j=1,\dots, d_i)$ and consider the
$0$--cycles $z(D^{i}_j)$. Then
$c_n(J_1(L))=\Sigma_{i=1}^s\Sigma_{j=1}^{d_i}z(D^{i}_j).$
\endproclaim

\demo{Proof} From \cite{LPS1, Cor. 2.6} we know that $c_n(J_1(L))$
is represented by $j_1(W)^{-1}(0)$.  For any $D^{i}_j \in |W| \cap
\Cal D^{i}$ and for any $x \in \text{Sing}(D^{i}_j)$ we can take
$s \in V$ the section defining $D^{i}_j$ and $t \in V$ not
vanishing at $x$ such that $W=<s,t>$. Now we can use the notation
of \cite{LPS2, Prop. 1.1}, i.e., there exist local coordinates
$x_1,\dots,x_n$ at $x \in X$ such that $s =\Sigma_{j=1}^n
a_{ij}x_ix_j+$h.o.t. (higher order terms) and
$j_1(t)=(1,0,\dots,0)$ since $t$ does not vanish at $x$. Hence
$j_1(t)\wedge j_1(s) =(\frac{\partial s}{\partial x_1}, \dots
,\frac{\partial s}{\partial x_n},0,\dots ,0).$ This shows that the
zero sub-scheme of $j_1(W)^{-1}(0)$ supported at $x$ is defined by
the Jacobian ideal $(\frac{\partial s}{\partial x_1}, \dots ,
\frac{\partial s}{\partial x_n}).$ Hence this 0--cycle has to be
$z_x=\mu_x(D^{i}_j)x$, where $\mu_x(D^{i}_j)$ is the Milnor number
of $x$ as an isolated singular point of $D^{i}_j$. This gives the
equality of 0--cycles of the statement. \qed \enddemo

From now on $c_n(J_1(L))$ will stand for the degree of the
corresponding 0--cycle. With the previous notation take a general
$D \in \Cal D^{l}$ where $1 \leq l\leq s$. Let us describe more
explicitly the singularities of $D$. By (1.3) (with the notation
there) and the classical biduality theorem, $D$ defines an element
of $\phi_V(Y_{i_lj_l})^\vee$ whose contact locus is a single point
$q \in \phi_V(Y_{i_lj_l})$. Then we claim that the singular locus
of $D$ is confined to $\phi_V^{-1}(q)=\{x_1,\dots,x_m\}$, that is,
$$\text{Sing}(D) \subseteq \phi_V^{-1}(q).  \tag 3.2.1$$
Since $D$ corresponds to an element of $\phi_V(Y_{i_lj_l})^\vee$
there exists $x \in \phi_V^{-1}(q)$ such that $D \in |V-2x|$.
Consider $y \neq x$ such that $D \in |V-2y|$. Since $D \in |V-2y|
\cap \Cal D^l$, $|V-2y| \subset \Cal D^l$. By (1.3) $y \in X_j$,
$j \geq i_l$. We have then an irreducible component $Y_{jk}
\subset {\overline {X}}_j$ such that $\dim(Y_{jk})=n-j$ and $\Cal
D^l=\phi_V(Y_{i_lk_l})^\vee=\phi_V(Y_{jk})^\vee$. If $j>i_l$, by
the classical biduality theorem, we get the  contradiction
$\phi_V(Y_{i_lk_l})=\phi_V(Y_{jk})$. Hence $y \in {\overline
X}_{i_l}$ and $\dim(|V-2y|)=N-1-(n-i_l)$. If $y$ is a smooth point
of $X_{i_l}$ then there exists $Y_{i_lk_l} \subset {\overline
{X}}_{i_l}$ with $\dim(Y_{i_lk_l})=n-i_l$ such that
$\phi_V(Y_{i_lj_l})^\vee= \phi_V(Y_{i_lk_l})^\vee$. Then
$\phi_V(Y_{i_lj_l})= \phi_V(Y_{i_lk_l})$ by the the classical
biduality theorem, as pointed out in (1.3.1). In particular
$\phi_V(x)=\phi_V(y)$ because the contact locus of a general
element on $\phi_V(Y_{i_lj_l})^\vee$ is just a point. This proves
the claim unless $y$ is a singular point of $\overline {X}_{i_l}$
and $\dim(|V-2y|)=N-1-(n-i_l)$. If this occurs, since
$\dim($Sing$(\overline{X}_{i_l}))<n-i_l$, we get the contradiction
$N-1=\dim(\Cal
D^l)=N-1-(n-i_l)+\dim($Sing$(\overline{X}_{i_l}))<N-1$.

By the previous discussion there exists an index $k$, $1 \leq k
\leq m$, such that, after reordering if necessary, $\text{\rm
Sing}(D)=\{x_1,\dots,x_k\} \subset X_{i_l} \cap \phi_V^{-1}(q)$
and $z(D)=\mu_{x_1}(D)x_1+\dots+\mu_{x_k}(D)x_k$. By
semicontinuity the degree of this 0--cycle is constant at the
general point of $\Cal D^l$. Let $m_l=\deg(z(D))$ for $D \in \Cal
D^l$ general. Hence, by (3.2) and recalling that
codeg$(X,V)=d_1+\dots+d_s$, we get
$$c_n(J_1(L))=m_1d_1+\dots+m_sd_s \geq \text{\rm codeg}(X,V). \tag
3.2.2$$ In particular, if the discriminant locus has just one
maximal dimensional irreducible component then there exists a
positive integer $m$ such that
$$c_n(J_1(L))= m\ \text{\rm codeg}(X,V). \tag 3.2.3$$

Classical results on Milnor numbers can be applied, for instance,
see \cite{DJP, Thm. 3.4.29, p. 122} :

\noindent (3.2.4) $z_x=x$ if and only if $x$ is an isolated
non--degenerate quadratic singularity of $D^{i}_{j}$.

In the classical case, see for example \cite{BS, Rmk. 1.6.11, p.
33}, when the dual variety is a hypersurface the general singular
hyperplane section has only an isolated non--degenerate quadratic
singularity, so that $c_n(J_1(L))=\text{\rm{codeg}}(X,V)$, being
$m=1$ in (3.2.3).

Let us give a bound when the dual of the image of $X$ by $\phi_V$
is the only maximal dimensional irreducible component of the
discriminant locus.

\proclaim{(3.2.5) Lemma} Let $(X,L,V)$ be as in $(0.0)$ with
$\dim(\Cal D)=N-1-k$. If $\phi_V(X)^\vee$ is the only
$(N-1-k)$--dimensional irreducible component of $\Cal D$ then
$c_n(J_1(L)) \leq \text{\rm
codeg}(X,V)\frac{L^n}{\deg(\phi_V(X))}. $
\endproclaim

\demo{Proof} If $k>0$ then $c_n(J_1(L))=0$ and the assertion is
obvious. If $k=0$ then (3.2.4) and the fact that the general
element of $\Cal D(X,V)$ has only isolated non--degenerate
quadratic singularities lead to the inequality $m \leq
\deg(\phi_V) =L^n/\text{deg}(\phi_V(X))$ in (3.2.3) and the
assertion follows. \qed
\enddemo

Let $(X,L,V)$ be as in $(0.0)$ with $\Cal J_2=\emptyset$ and
$\phi_V(\Cal J_1)$ not contained in the singular locus of
$\phi_V(X)$. As in \cite{BDL, Lemma 1 (3)}, for general $y \in
\phi_V(\Cal J_1)$ and $x \in \phi_V^{-1}(y)$ we can choose local
coordinates $x_1,\dots,x_n$ on $X$ centered at $x$ and
$y_1,\dots,y_n$ on $\phi_V(X)$ centered at $y$  to write
$y_1=x_1^{k},\; y_2=x_2,\; \dots,\; y_n=x_n,$ the branch locus
locally being defined by $y_1=0$. Then, an element $D \in \Cal
D^1$ singular at $x$ is defined by $s(x_1^k,\dots,x_n)$ where
$s(y_1,\dots,y_n)$ defines the corresponding hyperplane section
through $y$. Recall that $y=(0,\dots,0)$ is a smooth point of
$s(y_1,\dots,y_n)=0$ but $s(x_1^k,x_2,\dots,x_n)=0$ is singular at
$x=(0,\dots,0)$. Hence $s(y_1,\dots,y_n)=y_1+$h.o.t. In fact,
$s(x_1^k,x_2,\dots,x_n)=x_1^k+$h.o.t.

In particular, suppose $n=N=2$ and $\Cal J_2=\emptyset$. Since
$s(y_1,y_2)=0$ is the tangent line to the branch locus at the
(general) point $y=(0,0)$ then $s(y_1,y_2)=y_1+\alpha
y_2^2+$h.o.t., where $\alpha \in \Bbb C-\{0\}$. Hence:
$$s(x_1^k,x_2)=x_1^k+\alpha x_2^2+\text{h.o.t. and so } \mu_x(D)=k-1,
\tag 3.2.6$$ relating the index of ramification of $\phi_V$ at $x$
and the Milnor number of the singularity $x$. Of course $k \leq
{\text{\rm deg}}\phi_V \leq L^2$; then, specializing (3.2) in this
particular setting, we obtain the bound:
$$c_2(J_1(L))\leq {\text{\rm codeg}}(X,V)(L^2-1), \tag
3.2.7$$ which is sharp as example (b) in (3.3) will show.

If $\dim(\Cal D(X,V))=N-1-k$, $k>0$, then to get an expression as
in (3.2) we have to consider $z(D^{i}_j)$, the $k$--cycle
corresponding to the singular locus of a general point $D^{i}_j
\in \Cal D^{i}$, $\dim(\Cal D^i)=N-1-k$, counted with its
appropriate number. Generalizations of the Milnor number are
naturally considered, see \cite{A}. Let us produce some  examples
when $k=0$.

\noindent{\bf (3.3) Examples.}

{\bf (a)} Consider the example (1.1.0): $S$ is a a Del Pezzo
surface with $K_S^2=1$ and $L=-2K_S$. We know that $\phi_L:S \to
\Gamma \subset \Bbb P^3$ is a double cover of the quadric cone
$\Gamma$, branched at the vertex $v$ and along the smooth curve
$B$ cut out on $\Gamma$ by a transverse cubic surface. As
explained in (1.1.0) the maximal dimensional components of $\Cal
D(S,L)$ are $\Cal D_1=B^\vee$, of degree 18, and the plane $\Cal
D_2=v^\vee$. Note that $\Cal D_0=\gamma^\vee$ is a conic, the dual
of $\gamma$, contained in $\Cal D_1 \cap \Cal D_2$. Hence
$\text{codeg}(S,L) = \deg (\Cal D_1) + \deg (\Cal D_2) = 19$. On
the other hand, since $S$ is obtained by blowing-up $\Bbb P^2$ at
$8$ points, the Euler--Poincar\'e characteristic of $S$ is
$e(S)=11$. Moreover, $L^2=4$ and genus formula shows that
$g(L)=2$. Thus $c_2(J_1(L)) = 11 + 4 + 4 = 19$. Whence
$c_2(J_1(L))=\text{codeg}(X,V)$ which implies that the general
element in $\Cal D_1$ as well as that in $\Cal D_2$ has a single
isolated non-degenerate quadratic singularity.

{\bf (b)} Let $\pi:S \to \Bbb P^2$ be a cyclic cover of degree $d$
branched along a smooth curve $\Delta \in |\Cal O_{\Bbb
P^2}(bd)|$. Let $L:=\pi^* \Cal O_{\Bbb P^2}(1)$. We have $K_S =
\pi^* \Cal O_{\Bbb P^2}(b(d-1)-3)$ by the ramification formula,
and, for $i>0$ $h^i(\Cal O_S) = h^i(\Cal O_{\Bbb P^2}) + h^i(\Cal
O_{\Bbb P^2}(-b))+\dots+h^i(\Cal O_{\Bbb P^2}(-b(d-1)))$ by the
projection formula. Thus $q(S)=0$, while
$p_g(S)=\Sigma_{i=1}^{d-1}\binom{ib-1}{2}$. Therefore $12\chi(\Cal
O_S)=b^2(d-1)(d)(2d-1)-9bd(d-1)+12d.$ Since $K_S^2=
d(b(d-1)-3)^2$, Noether's formula and \cite{BS, Lemma 1.6.4} give
$$c_2(J_1(L))= 12 \chi(\Cal O_S) - K_S^2 +2K_S L+3L^2=
(d-1)bd(bd-1).$$ Since class$(\Delta)=$deg$(\Delta^\vee)=bd(bd-1)$
we have $c_2(J_1(L))=(d-1)$codeg$(S,L)$ and therefore the general
element of $\Cal D$ has only a single isolated non-degenerate
singularity whose Milnor number is $d-1$.

{\bf (c)} Let $(X,L,V)$ be as in (0.0) and suppose furthermore
that $L$ is very ample and $(X,L) \not= (\Bbb P^n, \Cal O_{\Bbb
P^n}(1))$. Let $h^0(L)=M+1$. Then $\Cal D(X,L) = \phi_L(X)^\vee
\subset \Bbb P^{M\vee}$. Suppose $\Cal D(X,L)$ is a hypersurface.
Then $c_n(J_1(L))$ represents its degree. Take $V$ general and let
$\dim(V) = N+1$. Note that $n \leq N \leq M$, the first inequality
following from the fact that $V$ spans $L$. Geometrically,
$\phi_V(X)$ is the general projection of $\phi_L(X) \subset \Bbb
P^M$ into $\Bbb P^N$ and in the dual projective space $\Cal
D(X,V)$ can be regarded as a general linear section of $\Cal
D(X,L)$. In fact $\Cal D(X,V) = \Cal D(X,L) \cap |V|$. Therefore
$\Cal D(X,V)$ is an irreducible (reduced) hypersurface of $|V|$,
since $V$ was chosen general. Moreover it has the same degree as
$\Cal D(X,L)$, hence $\deg(\Cal D(X,V))=c_n(J_1(L)),$ which
exactly means that the general element of $\Cal D(X,V)$ has an
isolated non-degenerate quadratic singularity. Moreover, if $n=N$
then $\Cal D^0=\emptyset$ and $\Cal D(X,V)$ is birational to $\Cal
J_1$. The birational map is given by $p_2 \circ p_1^{-1}$
recalling (0.1).

{\bf (d)} Let us consider some special projections of embedded
projective varieties. Let $X=\Bbb P^n$, $L=\Cal O_{\Bbb P^n}(m)$
and $V = < x_0^m, \dots , x_n^m >$. Then $V \subset H^0(X,L)$
spans $L$ and $\phi_V:X \to \Bbb P^n$ is a finite morphism of
degree $m^n$. Take a general pencil $\{D_t\}_{t \in \Bbb P^1}$ in
$|V|$. For every $t=(t_0:t_1) \in \Bbb P^1$, the hypersurface
$D_t$ has equation $a_0(t)x_0^m + \dots + a_n(t)x_n^m=0$, where
$a_i(t) = a_{i,0}t_0 + a_{i,1}t_1$ for $i=0,\cdots ,n$. Then $D_t$
is singular at $(x_0: \dots :x_n)$ if and only if
$a_0(t)x_0^{m-1}= \dots = a_n(t)x_n^{m-1}=0.$ We thus see, the
pencil being general, that there are exactly $n+1$ singular
hypersurfaces, say $D_0, \dots , D_n$, defined by $t$ satisfying
$a_i(t)=0$ for $i=0,\cdots ,n$ respectively. Each of them is
singular at one of the vertices of the homogeneous coordinate
system having a point of multiplicity $m$ there. E.\ g., $D_0$ has
equation $a_1(\tau)x_1^m + \dots + a_n(\tau)x_n^m=0$, where $\tau$
satisfies $a_0(\tau)=0$. In particular, this shows that
$\text{codeg}(X,V)=n+1.$ Moreover, $\Cal D(X,V)$ consists of $n+1$
hyperplanes. To see this, let $p$ be any of the points $(1:0:
\cdots :0), \dots, (0:0: \dots :1)$. Then each hyperplane $|V -
p|$ lies in $\Cal D(X,V)$ (and in fact $|V - p| = |V - mp|$ for
each of them). Actually, let $a=(a_0: \dots :a_n)$; then the
hypersurface $D_a$ has equation: $a_0x_0^m + a_1x_1^m + \dots +
a_nx_n^m = 0$. Let $p = (1:0: \dots :0)$. Then $D_a \in |V - p|$
if and only if $a_0x_0^m=0$, i.e., if $a_0=0$. Thus $|V - p|$
consists exactly of the hypersurfaces of equation $a_1x_1^m +
\dots + a_nx_n^m=0$, which have $p$ as a singular point of
multiplicity $m$. Let us prove also the following

\proclaim{(3.3.1) Proposition} For $(\Bbb P^n,L=\Cal O_{\Bbb
P^n}(m),< x_0^m, \cdots , x_n^m >)$, $c_n(J_1(L)) = (m-1)^n
\text{\rm{codeg}}(X,V).$
\endproclaim

\demo{Proof} Let $p=(1:0: \dots :0)$ and $D \in |V-p|$. Let
$y_1,\dots,y_n$ be local coordinates of $X$ at $p$. Arguing as
after (3.1) and noting that $j_1(s) =(s, ma_1 y_1^{m-1}, \dots
,ma_n y_n^{m-1})$, where $s$ defines $D$, we get
$\mu_p(D)=(m-1)^n$. Since the same computation can be done for any
point of the set $\{(1:0:\dots:0),\dots, (0:\dots:1)\}$ we get the
final assertion. \qed \enddemo

The previous computation can be done in a different way. Let $X$
be any projective manifold of dimension $n$, let $L$ be a line
bundle on $X$ and consider the exact sequence: $0 \to \Omega_X^1
\otimes L \to J_1(L) \to L \to 0.$ We have (e.\ g., see \cite{BS,
Lemma 1.6.4}) $c_n(J_1(L))= \sum_{i=0}^n
\binom{n+1-i}{n-i}c_i(\Omega_X^1) L^{n-i} = \sum_{i=0}^n (n+1-i)
c_i(\Omega_X^1) L^{n-i}.$ So, for $X=\Bbb P^n$, letting $h = \Cal
O_{\Bbb P^n}(1)$ and recalling that the total Chern class of
$\Omega_X^1$ is $(1-h)^{n+1}$ mod $h^{n+1}$, we get $c_n(J_1(L))=
\sum_{i=0}^n (n+1-i) (-1)^i \binom{n+1}{i} h^i L^{n-i}.$ An
immediate check shows that $(n+1-i) \binom{n+1}{i} = (n+1)
\binom{n}{i}.$ Hence, for $L= \Cal O_{\Bbb P^n}(m)$, we get
$$c_n(J_1(L)) = (n+1) \sum_{i=0}^n \binom{n}{i} m^{n-i} (-1)^i =
(n+1) (m-1)^n. \tag 3.3.2$$

This example is nice also from the point of view of the biduality
theorem. Let $p_i=(0: \dots :0:1:0 \dots:0)$ be the
$i^{\text{th}}$ vertex of the coordinate system, $i=0,\dots,n$.
The first jumping set $\Cal J_1$ is the union of the coordinate
hyperplanes $x_i=0$, and it is easy to realize that each of the
further jumping sets is the union of all coordinate linear
subspaces of the appropriate dimension; in particular, $\Cal J_n$
consists of the $n+1$ points $p_0, \dots , p_n$. Now, let $D \in
\Cal D(X,V)$ be defined by the section $s= \sum_{i=0}^n a_i x_i^m
\in V$ and suppose that the singular locus of $D$ includes a point
$p$, distinct from the $p_i$'s. Then, up to reordering the
coordinates, $p=(0: \dots :0:y_s: \dots :y_n)$, with $y_j \not= 0$
for $j \geq s$, and
$$a_{s}= \dots = a_n = 0. \tag 3.3.3$$
Therefore $D$ lies on the intersection of $n-s+1$ of the $n+1$
hyperplanes of $\Bbb P^{n \vee}$ constituting $\Cal D(X,V)$. On
the other hand, according to our choice, $p$ lies on an
irreducible component $Y$ (a linear space) of the jumping set
$\Cal J_s$ (since it is required that $s$ coordinates vanish).
Note that here $\Cal J_s$ is the same as $X_s$ (the closure of
$\Cal J_s \setminus \Cal J_{s+1}$, as defined in Section 0).
Moreover, $\phi_V(Y)^{\vee}$ is the $\Bbb P^{s-1} \subset |V|$
defined by (3.3.3), which is contained in $\Cal D_s(X,V) \subset
\Cal D(X,V)$. This discussion illustrates (1.2) very well.
Moreover, letting $s=n$ we get a significant example also for
(1.3): the index named $i$ is $n$ and the corresponding $Y$ is
simply the point $p_n$.

Let us now focus on low codegree triplets. By \cite{LPS1, Th. 2.8}
the only $(X,L,V)$ as in (0.0) such that codeg$(X,V)=0$, which
means $\Cal D(X,V)=\emptyset$, is $(\Bbb P^n, \Cal O_{\Bbb
P^n}(1),H^0(\Bbb P^n, \Cal O_{\Bbb P^n}(1)))$. It becomes then
natural, as in the classical case, to face the problem of
classifying low codegree triplets $(X,L,V)$. Let us first show
some examples.

\flushpar{\bf (3.4) Examples.}

{\bf (a)} When considering smooth projective varieties $X
\subseteq \Bbb P^N$, the only codegree 1 varieties are
(degenerate) linear spaces $\Bbb P^n \subset \Bbb P^N$ with $n<N$,
and the only codegree 2 varieties are quadrics. That is, there are
no triplets $(X,L,V)$ as in (0.0) with $\phi_V$ an embedding and
$\text{codeg}(X,V)=1$ and the only example with
$\text{codeg}(X,V)=2$ is $(Q,\Cal O_Q(1), H^0(Q,\Cal O_Q(1)))$
where $Q \subset \Bbb P^{n+1}$ is a smooth quadric. If
codeg$(X,V)=3$ and $\phi_V$ is an embedding a complete
classification can be found in \cite{Z1}, \cite{Z2, IV.5}.

{\bf (b)} Take a cyclic covering $f:X \to \Bbb P^n$, branched
along a smooth quadric $Q \subset \Bbb P^n$. For the triplet $(X,
L=f^*\Cal O_{\Bbb P^n}(1), V=f^*H^0(\Bbb P^n,\Cal O_{\Bbb
P^n}(1)))$ the discriminant locus is a smooth quadric, $\Cal
D(X,V)=Q^\vee$. In fact, $(X, L)$ is as in (a) above, but $V$ is a
codimension 1 general linear subspace of $H^0(X,L)$.

{\bf (c)} For the example (b.1) of (1.1.1) $\Cal D(X,V)=\Bbb P^1
\times \Bbb P^{n-1} \subset \Bbb P^{2n-1}$ and so codeg$(X,V)=n$.

{\bf (d)} Consider now $(\Bbb P^2, \Cal O_{\Bbb P^2}(2))$. If
$V=H^0(\Bbb P^2, \Cal O_{\Bbb P^2}(2))$ then $\Cal D(\Bbb P^2,V)$
is the dual variety of the Veronese surface $S \subset \Bbb P^5$,
that is the cubic symmetroid $S^\vee \subset \Bbb P^5$, hence the
codegree is 3. In fact, taking coordinates $x_0,\dots,x_5$ in
$\Bbb P^5$, $S^{\vee}$ is defined by $rk(M) \leq 2$ being
$M=\pmatrix
x_0 &x_3 &x_4 \\
x_3 & x_1 & x_5\\
x_4 &x_5 &x_2\\ \endpmatrix$ and its singular locus is a Veronese
surface defined by $rk(M)=1$.

Let $V \subset H^0(\Bbb P^2, \Cal O_{\Bbb P^2}(2))$ be a subspace
such that $|V|$ is a base-point free web of conics and consider
the morphism $\phi_V: \Bbb P^2 \to \Bbb P^3$. Since
$\Sigma:=\phi_V(\Bbb P^2)$ is non-degenerate, there are only two
possibilities: either $\Sigma$ is a quartic surface and $\phi_V$
is birational, or $\Sigma$ is a quadric surface and $\phi_V$ has
degree $2$. Note that in the latter case  $\Sigma$ must be a
quadric cone. Actually, $\Bbb P^2$ cannot be a branched double
cover of a smooth quadric surface; otherwise, the ramification
formula would imply that $9=K_{\Bbb P^2}^2 \in 2\Bbb Z$, a
contradiction. Here are examples of both cases:
\par (d.1) Let $(u:x:y)$ be homogeneous coordinates in
$\Bbb P^2$ and let $V_1 = \langle u^2+x^2+y^2, xy, uy, ux \rangle$
Then $|V_1|$ is a base-point free web of conics and
$\phi_{V_1}:\Bbb P^2 \to \Bbb P^3$ is a birational morphism onto
the roman Steiner quartic surface $\Sigma$ of equation
$y_1^2y_2^2+y_1^2y_3^2+y_2^2y_3^2- y_0y_1y_2y_3=0$, where
$(y_0:y_1:y_2:y_3)$ are the corresponding homogeneous coordinates
in $\Bbb P^3$. Note that $\Cal D(\Bbb P^2, V_1)$ is the
irreducible cubic surface defined by: $\lambda^3 - \lambda(\mu^2 +
\nu^2 + \epsilon^2) +2 \mu \nu \epsilon = 0$, where $(\lambda:
\mu: \nu:\epsilon)$ are the dual homogenous coordinates in $\Bbb
P^{3 \vee}$.

\par  (d.2) Let $V_2=\langle u^2, x^2, xy, y^2\rangle$. Then
$|V_2|$ is a base-point free web of conics and the morphism
$\phi_{V_2}:\Bbb P^2 \to \Bbb P^3$ is two--to--one onto the
quadric cone $\Sigma$ of equation $y_1y_3-y_2^2=0$.  Note that
$\Cal D(\Bbb P^2, V_2)$ is defined by the equation $\lambda(\mu
\epsilon - \nu^2)=0$, hence it is reducible into a plane plus a
quadric cone.

Now let us choose $V \subset H^0(\Bbb P^2, \Cal O_{\Bbb P^2}(2)))$
such that $|V|$ is a base-point free net of conics. By \cite{W}:

(d.3) either $\Cal D(\Bbb P^2, V) \subset \Bbb P^2$ is irreducible
and there exists a suitable choice of homogeneous coordinates
$x,y,z$ in the plane (suggested in \cite{W}) so that
$$V=\langle 2xz+y^2,2yz,-x^2-2gy^2+cz^2+2gxz
\rangle \subset H^0(\Bbb P^2, \Cal O_{\Bbb P^2}(2)), \tag 3.4.1$$
where $g,c$ are complex parameters, or

(d.4) $\Cal D(\Bbb P^2, V) \subset \Bbb P^2$ is reducible and for
a suitable choice of homogeneous coordinates, either $V=\langle
x^2,y^2,z^2 \rangle$ or $V=\langle x^2,y^2, z^2+2xy \rangle$.

Then, the possibilities are the following:

(d.3.1) If $c \ne -9g^2$ in (3.4.1) then $\phi_V: \Bbb P^2 \to
\Bbb P^2$ is a degree four map branched along a sextic $C \subset
\Bbb P^2$ with 9 cusps (and no other singularities). Hence $\Cal
D(\Bbb P^2, V)=\Cal D_1(\Bbb P^2,V)=C^\vee \subset \Bbb P^2$ is a
smooth plane cubic.

(d.3.2) If $c=-9g^2 \ne 0$ in (3.4.1) then $\phi_V: \Bbb P^2 \to
\Bbb P^2$ is a degree four map branched along a quartic curve $C
\subset \Bbb P^2$ with three cusps (and no other singularities).
Hence $\Cal D(\Bbb P^2,V)=\Cal D_1(\Bbb P^2,V)=C^\vee$ is a nodal
plane cubic.

Note that for $c=g=0$ $|V|$ is not base point free,
Bs$|V|=\{(0:0:1)\}$.

(d.4.1) If $V=\langle x^2,y^2, z^2+2xy \rangle$ then $\phi_V$ is a
degree four map branched along the union of a smooth conic $Q$ and
two lines tangent to $Q$. Hence $\Cal D_1(\Bbb P^2,V)=Q^\vee$,
$\Cal D_2(\Bbb P^2,V)=\ell$ and $\Cal D(\Bbb P^2,V)=\Cal D_1(\Bbb
P^2,V) \cup \Cal D_2(\Bbb P^2,V)$ where $Q^\vee$ is a smooth conic
and $\ell$ is a line transverse to $Q^\vee$.

(d.4.2) If $V=\langle x^2,y^2,z^2 \rangle$ then $\phi_V:\Bbb P^2
\to \Bbb P^2$ is a degree four map branched along three general
lines. Hence $\Cal D(\Bbb P^2,V)=\ell_1 \cup \ell_2 \cup \ell_3
\subset \Bbb P^2$, being $\ell_1,\ell_2$ and $\ell_3$ the dual
trilateral to the branch locus. Hence $\Cal D(\Bbb P^2,V)=\Cal
D_2(\Bbb P^2,V)=\{\ell_1,\ell_2,\ell_3\}$ and $\Cal D_1(\Bbb
P^2,V)=\{(\ell_1 \cap \ell_2), (\ell_2 \cap \ell_3), (\ell_1 \cap
\ell_3)\}$.

A corollary of the study of nets of conics according to \cite{W}
is the following. For a suitable $3$-dimensional vector subspace
$V \subset H^0(\Bbb P^2, \Cal O_{\Bbb P^2}(2))$ it may happen that
the corresponding plane section of the cubic symmetroid either is
a cuspidal curve, or contains a double line, or is the union of a
smooth conic and one of its tangent lines. In all these cases,
however, $|V|$ is not base-point free.

We can exclude the possibility for the codegree to be one.

\proclaim{(3.5) Theorem} Let $(X,L,V)$ be as in $(0.0)$. Then {\rm
codeg}$(X,V) \geq 2$ except for $(X,L,V)=(\Bbb P^n,\Cal O_{\Bbb
P^n}(1),H^0(\Bbb P^n,\Cal O_{\Bbb P^n}(1))$ (for which {\rm
codeg}$(X,V) =0$).
\endproclaim

\demo{Proof} By assumption {\rm codeg}$(X,V)=0$ if and only if
$\Cal D(X,V)=\emptyset$, which implies $\dim(\Cal D(X,V))=-1$. We
conclude by \cite{LPS, Thm. 2.8}. Let codeg$(X,V)=1$ then $\Cal
D(X,V)=\Cal T \cup \Cal D'$ where $\Cal T =\Bbb P^{N-1-k}$ and
$\Cal D'$ is the union of the irreducible components of the
discriminant locus of dimension $<N-1-k$. By (1.4.1) there exists
$Y \subseteq \Cal J_{n-k}(X,V)$ such that $\phi_V(Y)=\Bbb P^{k}$,
$\Cal T=\phi_V(Y)^\vee$ and $\Cal J_{n-k+i}=\emptyset$ for $i>0$.
If $\Cal D' \ne \emptyset$ consider $\Cal D^0$ the union of the
maximal dimensional components of $\Cal D'$. By (1.3) $\Cal
D^0=\phi_V(Y_1)^\vee \cup \dots \cup \phi_V(Y_s)^\vee$. Suppose
dim$(\Cal D^0)=N-1-k_0$. Take general sections $s_1, \dots,
s_{k_0+1} \in V$ and let $M=|\langle s_1,\dots,s_{k_0+1}\rangle|
\subseteq |V|$. Since $\dim(M \cap \Cal D^0)=0$ the classical
biduality theorem implies that the image by $\phi_V$ of the
singular locus of any element corresponding to a point $p \in M
\cap \Cal D^0$ is a linear space of dimension $k_0$, say $\Bbb
P^{k_0}_p$. Then, by \cite{LPS1, Thm. 2.4}, $\phi_V((j_1(s_1)
\wedge \dots \wedge j_1(s_{k_0+1}))^{-1}(0))$ is a finite union of
linear spaces of dimension ${k_0}$.  By \cite{LPS1, 2.3.2}
$\phi_V(\Cal J_{n-k_0})$ is the intersection of these linear
spaces. Since $\phi_V(Y) = \Bbb P^k \subseteq \phi_V(\Cal J_{n-k})
\subseteq \phi_V(\Cal J_{n-k_0})$ it follows that all the contact
loci of $\phi_V(Y_i)$ ($1 \leq i \leq s$) are meeting along
$\phi_V(Y)$. This in particular means that the intersection of all
the projective tangent spaces to $\phi_V(Y_i) \subset \Bbb P^N$ is
not empty. So, $\phi_V(Y_i)$ is a cone \cite{R, Prop. 1.2.6} whose
vertex contains $\phi_V(Y)$, contradicting the non-emptyness of
$\Cal D'$. Hence $\Cal D'=\emptyset$ so that $\Cal D$ is
irreducible. This implies $\phi_V(X)$ is a cone whose vertex
contains $\phi_V(Y)$ and $\phi_V(\Cal J_1)$ is a union of cones
with vertex containing $\phi_V(Y)$ as in (2.7.1). In this
situation the general line in $\Bbb P^n$ only meets $\phi_V(\Cal
J_1)$ at one point. This leads to a contradiction as in the proof
of (2.7). \qed
\enddemo

\head 4 Low codegree curves
\endhead

In this section we classify curves of codegree less than or equal
to three.

\proclaim{(4.1) Remark} For $(X,L,V)$ as in $(0.0)$
$\phi_V(X)^\vee$ cannot be a  cone.
\endproclaim

In fact, if $\phi_V(X)^\vee$ is a cone, then $\phi_V(X)^{\vee
\vee}=\phi_V(X)$ is degenerate, contradicting the assumptions of
(0.0). We collect two basic facts on curves in the following
remark. Proofs are straightforward.

\proclaim{(4.2) Remark} Let $(C,L,V)$ be as in $(0.0)$ with
$\dim(C)=1$. Then:

$(4.2.1)$ $\phi_V(C)^\vee$ is either empty or a hypersurface of
$\Bbb P^N$ of degree $\geq 2$,

$(4.2.2)$ $\Cal D(C,V)$ contains $|\phi_V(\Cal J_1)|$ hyperplanes.

\endproclaim

If codeg$(C,V)=1$ then $\phi_V(C)=\Bbb P^1$ and $|\phi_V(\Cal
J_1)|=1$ by (4.2), contradicting (2.3). This gives a different
proof of (3.5) in the case of curves.

If codeg$(C,V)=2$ then either $\phi_V(C)^\vee$ is a quadric and
$\Cal J_1=\emptyset$, or $\phi_V(C)=\Bbb P^1$ and $\phi_V(\Cal
J_1)=\{p_1,p_2\}$. If the former holds then $\phi_V(C)^\vee$ is
smooth by (4.1), and so, by biduality, $\phi_V(C)$ is a smooth
conic. Since $\Cal J_1=\emptyset$ then $(C,L,V)=(\Bbb P^1, \Cal
O_{\Bbb P^1}(2), H^0(\Bbb P^1, \Cal O_{\Bbb P^1}(2)))$. If the
latter holds then, arguing as in (2.3), from the Riemann--Hurwitz
formula we obtain $2g(C)-2=-m_1-m_2$ being
$m_1=|\phi_V^{-1}(p_1)|$ (respectively $m_2=|\phi_V^{-1}(p_2)|$).
Then $g(C)=0$ and $m_1=m_2=1$. In fact, there exists an integer
$r>1$ such that $C=\Bbb P^1$, $L=\Cal O_{\Bbb P^1}(r)$ and $V$ has
to be chosen in the following way: since the complete linear
system embeds $\Bbb P^1$ in $\Bbb P^r$ as a rational normal curve,
we have to project it from $T=\Bbb P^{r-2}$ in such a way that
$|\phi_V(\Cal J_1)|=2$. Then $T$ is the intersection of two linear
spaces of dimension $r-1$ that are $(r-1)$-osculating to the
rational normal curve. This concludes the codegree 2 case.

If codeg$(C,V)=3$ then either

(4.3.1) $\phi_V(C)^\vee$ is a cubic and $\Cal J_1=\emptyset$, that
is, $\phi_V$ is an immersion, or

(4.3.2) $\phi_V(C)^\vee$ is a conic and exists $p \in \phi_V(C)$
such that $\phi_V(\Cal J_1)=p$, or

(4.3.3) $\phi_V(C)=\Bbb P^1$ and there exist three distinct points
$p_1,p_2,p_3 \in \phi_V(C)$ such that $\phi_V(\Cal J_1
)=\{p_1,p_2,p_3\}$.

\noindent If (4.3.1) holds, $\phi_V$ is an immersion and then
$\Cal D(X,V)=\phi_V(C)^\vee$. By (3.2.5) $c_1(J_1(L)) \leq
3\text{deg}(\phi_V)$. Since $\phi_V(C)^{\vee}$ is a cubic,
$\text{deg}(\phi_V(C)) \geq 3$. The previous bounds and the fact
that
$c_1(J_1(L))=2g(C)-2+2\text{deg}(L)=
2g(C)-2+2\text{deg}(\phi_V)\text{deg}(\phi_V(C))$
leads to a contradiction. Hence this case does not occur. If
(4.3.2) holds then $2g(C)-2=-d-m$, where $m=|\phi_V^{-1}(p)|$, by
the Riemann--Hurwitz formula. This gives $g(C)=0$ and $d=m=1$, a
contradiction. If (4.3.3) holds then, just as before,
Riemann--Hurwitz formula says $2g(C)-2=d-(m_1-m_2-m_3)$ where
$m_i=|\phi_V^{-1}(p_i)|$ for $i=1,2,3$. Whence:

\proclaim{(4.4) Theorem} Let $(X,L,V)$ as in $(0.0)$ with
$\dim(X)=1$ and \text{\rm codeg}$(X,V)\leq 3$. Then, either

$(4.4.1)$ $(X,L,V)=(\Bbb P^1, \Cal O_{\Bbb P^1}(2),H^0(\Bbb P^1,
\Cal O_{\Bbb P^1}(2)))$, \text{\rm codeg}$(X,V)=2$, or

$(4.4.2)$ $(X,L,V)=(\Bbb P^1, \Cal O_{\Bbb P^1}(r),V)$ with $r
\geq 2$ and $V \subset H^0(\Bbb P^1, \Cal O_{\Bbb P^1}(r))$ is
such that $\phi_V$ is the projection of $\phi_L(\Bbb P^1) \subset
\Bbb P^r$ from the intersection of two $(r-1)$--dimensional linear
spaces of $\Bbb P^r$ that are $(r-1)$--osculating to $\phi_L(\Bbb
P^1)$; \text{\rm codeg}$(X,V)=2$, or

$(4.4.3)$ $\phi_V(X)=\Bbb P^1$, $\phi_V(\Cal J_1)=\{p_1,p_2,p_3\}$
and \text{\rm codeg}$(X,V)=3$.
\endproclaim

\noindent{\bf (4.5) Examples} (showing that the list of (4.4) is
effective).

{\bf (a)} Consider a degree $r$ rational normal curve $C \subset
\Bbb P^r$ and for every $k \leq r-1$ let $\text{Osc}_c^k(C)$ be
the $k$-th osculating space to $C$ at $p$. Take two general points
$p_1$, $p_2 \in C$ and consider $M=\text{Osc}_{p_1}^{r-1}(C)\cap
\text{Osc}_{p_2}^{r-2}(C)=\Bbb P^{r-3}$. Let $T$ be a general
$\Bbb P^{r-2}$ in $\text{Osc}_{p_1}^{r-1}(C)$ containing $M$. Then
the projection from $T$ ramifies in $p_1$ with ramification index
$r$ and in $p_2$ with ramification index $r-1$. Since
$-2=-2r+(r-1+r-2+1)$ then the projection from $T$ ramifies in a
third point $p_3$ with ramification index 2 and we are in case
(4.4.3).

{\bf (b)} With the notation of (4.3.3) let us construct an example
with $g(C)=1$, $d=3$ and $m_1=m_2=m_3=1$. Take the projection of a
smooth plane cubic $C \subset \Bbb P^2$ from a point $x \in \Bbb
P^2 \setminus C$ onto $\Bbb P^1$. In order to have codegree three
we have to choose $x$ in the intersection of three tangent lines
to $C$ at flexes of $C$. For example consider the cubic defined by
the equation $x_0^3-x_1x_2^2+x_1^2x_2=0$ and project from
$(1:0:0)$ which is on the intersection of the tangent lines to the
cubic at the three flexes $(0:0:1)$, $(0:1:0)$ and $(0:1:1)$.

\head 5 Low codegree surfaces: general facts
\endhead

Consider now a triplet $(S,L,V)$ as in (0.0), $S$ being a surface.
Recall  that $\dim(\Cal D(S,V))=N-1$, see \cite{LPS1, Thm.\ 2.8}.
Suppose that $\phi_V(S)^\vee$ is the only $(N-1)$--dimensional
irreducible component of $\Cal D(X,V)$. If $\text{codeg}(S,V) <
\deg(\phi_V(S))$ then (3.2.5) combined with \cite{LPS1, Prop. A.1}
gives  the bound $L^2 - 1 \leq c_2(J_1(L)) \leq \text{codeg}(S,V)
 \frac{L^2}{\deg(\phi_V(S))} < L^2.$ It follows that the
first inequality has to be an equality, and so \cite{LPS1, Prop.\
A.1} implies that $(S,L) = (\Bbb P^2, \Cal O_{\Bbb P^2}(2))$. In
particular, $(S,V)$ is one of the pairs discussed in (d) of (3.4).
It turns out that, apart from this case, the inequality
$\text{codeg}(S,V) < \deg(\phi_V(S))$ cannot be true. Moreover, if
equality holds then $(S,L)$ is a scroll by \cite{LPS1, Prop.\ A1}.
So this proves the following

\proclaim{(5.1) Corollary} Let $(S,L,V)$ be as in $(0.0)$, where
$\dim(S)=2$. If $\phi_V(S)^\vee$ is the only $(N-1)$--dimensional
irreducible component of $\Cal D(X,V)$ then, either

$(5.1.1)$ $(S,L) = (\Bbb P^2, \Cal O_{\Bbb P^2}(2))$, or

$(5.1.2)$ $\text{\rm{codeg}}(S,V) \geq \deg(\phi_V(S))$.

\noindent Moreover if equality holds in $(5.1.2)$ and we are not
in $(5.1.1)$ then $(S,L)$ is a scroll.
\endproclaim

We can regard (5.1) as a natural extension of classical results of
Marchionna \cite{M} and Gallarati \cite{G} to the ample and
spanned setting. Let us recall here that for a triplet $(S,L,V)$
as in (0.0) with $\dim(S)=2$, $L$ very ample and $\phi_V$ an
embedding, it is usual to use the term {\it class} to refer to
codeg$(S,V)$. Marchionna proved that the class of a surface is
greater than or equal to its degree minus one and equality holds
when $(S,L)=(\Bbb P^2, \Cal O_{\Bbb P^2}(t))$, $t=1,2$, see
\cite{M}. Moreover Gallarati showed that the class is equal to the
degree if and only if $(S,L)$ is a scroll, see \cite{G}. The
example (d) of (2.6) is interesting in connection with (5.1). In
fact, for the elliptic scroll of invariant $e=-1$ when considering
$L=C_0+f$, $V=H^0(S,L)$ we obtain codeg$(S,V)=3$ and
$\phi_V(S)=\Bbb P^2$. On the other hand, when considering
$L=C_0+2f$, $V=H^0(S,L)$ then $\phi_V$ is an embedding and
$\phi_V(S) \subset \Bbb P^4$ is the quintic elliptic scroll.
Whence codeg$(S,V)=5$.

As said in the introduction, the geometry of $\phi_V(S) \subseteq
\Bbb P^N$ is an important tool in the study of $\Cal D(S,V)$. In
particular $\phi_V(S)^\vee \subseteq \Cal D(X,L)$ is a relevant
part of the discriminant. Let us comment some consequences of
$\phi_V(S)^\vee$ to be small. We will recall the following
definition and notation:

\proclaim{Definition} Let $Y \subset \Bbb P^M$ a projective
variety. The {\rm tangent developable to $Y$} is denoted by $TY$
an is defined as the closure in $\Bbb P^M$ of the union of the
embedded projective tangent spaces to $Y$ at its smooth points.
\endproclaim

\proclaim{(5.2) Proposition} Let $(S,L,V)$ be a triplet as in
$(0.0)$ with $\dim(S)=2$ and $\dim(\phi_V(S)^\vee)<N-1$. Then:

$(5.2.1)$ either $\phi_V(S) \subseteq \Bbb P^N$ is a cone, or

$(5.2.2)$ there exists a curve $C \subset \Cal J_1$ such that
$\text{\rm deg}(\phi_V(S)) \leq \text{\rm deg}(\phi_V(C)^\vee)\leq
\text{\rm codeg}(S,V)$ and $\phi_V(S)=T\phi_V(C)$.

\noindent In particular if $(5.2.1)$ does not hold then {\rm
codeg}$(S,V)>3$.
\endproclaim

\demo{Proof} Suppose $\phi_V(S) \ne \Bbb P^2$, if not we are in
case (5.2.1). Since the dual of $\phi_V(S) \subset \Bbb P^N$ is
not a hypersurface then the general tangent hyperplane is tangent
to $\phi_V(S)$ along a line. In particular $\phi_V(S) \subset \Bbb
P^N$ is swept out by lines. Since $\phi_V(S) \ne \Bbb P^2$ there
is a finite number of lines through the general point of
$\phi_V(S)$. Consider $x \in \phi_V(S)$ general. The general
tangent hyperplane $H$ to $\phi_V(S)$ at $x$ is tangent along a
line $\ell_H$ through $x$. Since there is a finite set of lines on
$\phi_V(S)$ through $x$ it holds that $\ell_H=\ell$ for the
general $H$ containing $T_{\phi_V(S),x}$. This says that
$T_{\phi_V(S),x}=T_{\phi_V(S),y}$ for a general $y \in \ell$. In
particular $\phi_V(S) \subset \Bbb P^N$ is a developable surface.
Then, by \cite{FP, Thm. 2.2.8}, either $\phi_V(S) \subset \Bbb
P^N$ is a cone, and so (5.2.1) holds, or it is the tangent
developable to a curve $E \subset \phi_V(S)$. Suppose
$\phi_V(S)=TE$. Then the general line in $\phi_V(S)$ is the
tangent line to $E$ at a smooth point $e \in E$, say $T_{E,e}$.
The general hyperplane $H$ containing this line cuts out
$\phi_V(S)$ along a reducible curve by degree reasons, and so its
corresponding element $D \in |V|$ is reducible, hence singular by
(0.4). This implies $E^{\vee} \subset \Cal D(S,V)$. In particular,
it is a non-linear $(N-1)$--dimensional irreducible component of
$\Cal D(S,V)$ and we conclude the existence of a curve $C \subset
\Cal J_1$ by (1.3).

Consider now a general line $R$ contained in $|V|$ corresponding
to the hyperplanes in $\Bbb P^N$ containing a fixed $T=\Bbb
P^{N-2}$. Since $E^{\vee}$ is an irreducible component of the
discriminant then there exist points $e_1,\dots,e_r \in E$,
$r=\text{deg}(E^\vee)\leq \text{codeg}(S,V)$, such that $\dim
(\langle T_{E,e_i}, T\rangle)=N-1$ for $1 \leq i \leq r$. This is
equivalent to saying that $T \cap T_{E,e_i}=x_i \ne \emptyset$ for
$1 \leq i \leq r$. In particular $\{x_1,\dots,x_r\} \subseteq T
\cap \phi_V(S)$. We claim that this is an equality. Indeed, if
there exists $x \in (T \cap \phi_V(S)) \setminus
\{x_1,\dots,x_r\}$ it holds that there exists an analytic arc
$\{e(t)\} \subset E$ such that $x \in \ell \subset \phi_V(S)$,
being $\ell$ the limit of the tangent lines $T_{E,e(t)}$. Let us
remark that if $\ell = T_{E,e_i}$ for some $i$ and $x \in \ell =
T_{E,e_i} \setminus\{x_1,\dots,x_r\}$ then $|\ell \cap T| \geq 2$
and so $\ell \subset T \cap \phi_V(S)$, contradicting the general
choice of $T$. Then $\ell \cap T =\{x\} \ne \emptyset$ and the set
of hyperplanes containing $\ell$ is contained in $E^\vee$, a
contradiction with deg$(E^\vee)=r$. Since developable quadrics and
cubics have to be cones, see \cite{E, pp.\ 32--33}, then
codeg$(X,V)>3$.\qed
\enddemo

For low codegree we can study the possibility for any maximal
dimensional component of the discriminant to be linear.

\proclaim{(5.3) Lemma} Let $(S,L,V)$ be as in $(0.0)$ and
$\dim(S)=2$. If any maximal dimensional component of $\Cal D$ is
linear then  $\phi_V(S)=\Bbb P^2$, $\phi_V(\Cal J_1)$ is a union
of at least three not collinear lines and $\phi_V(p)$ is contained
in a line of $\phi_V(\Cal J_1)$ for any $p \in \Cal J_2$. Moreover
$\text{\rm codeg}(S,V) \geq 3$ and if equality holds then
$(S,L)=(\Bbb P^2, \Cal O_{\Bbb P^2}(r))$.
\endproclaim

\demo{Proof} Since $\phi_V(S)^\vee \subseteq \Cal D(S,V)$, (5.2)
applies. Then either $\phi_V(S) \subseteq \Bbb P^N$ is a cone or
$\phi_V(S)$ is a developable surface (different from a cone). In
the second case by (5.2.2) there exists a curve $C \subseteq \Cal
J_1$ such that $T\phi_V(C)=\phi_V(S)$, in particular
$\phi_V(C)\subset \Bbb P^N$ is non-degenerate, and
$\phi_V(C)^\vee$ is a component of $\Cal D(S,V)$. This is a
contradiction because the dual of any non-degenerate curve is a
non-linear hypersurface. Hence $\phi_V(S) \subseteq \Bbb P^N$ is a
cone.

If $\phi_V(S) \subseteq \Bbb P^N$ is a not linear cone then the
vertex is a point, say $v$. If $\phi_V(\Cal J_1)$ contains an
irreducible curve of degree $\geq 2$ then its dual is a nonlinear
component of dimension $N-1$ of $\Cal D$, a contradiction. Then
$\phi_V(\Cal J_1)$ is a union of lines through $v$. This gives
$(S,L)=(\Bbb P^2,\Cal O_{\Bbb P^2}(1))$ by exactly the same
argument as after (2.7.1). In fact, the preimage of a general line
$\ell \subset \phi_V(S)$ through $v$ is a curve whose singular
locus is contained in $\phi_V^{-1}(v)$. Then for the normalization
$\gamma$ of any of its irreducible components, $\phi_V$ defines a
map from $\gamma$ onto $\Bbb P^1$ branched only at $v$. By (2.3)
any irreducible component of $\phi_V^{-1}(\ell)$ is then
isomorphic to $\Bbb P^1$ via $\phi_V$. Hence $(S,L)$ is swept out
by lines and we conclude by \cite{LP2, Thm. 1.4}. Hence
$\phi_V(S)=\Bbb P^2$.

Since $\phi_V(S)=\Bbb P^2$, $\Cal J_1$ is a union of curves and so
$\phi_V(\Cal J_1)$ is a union of lines, being linear any maximal
dimensional  component of the discriminant. Moreover, since $\Cal
J_2 \subset \Cal J_1$, any $p \in \Cal J_2$ is contained in a
curve of $\Cal J_1$, hence $\phi_V(p)$ is contained in a line of
$\phi_V(\Cal J_1)$. By exactly the same argument of the previous
paragraph, $\phi_V(\Cal J_1)$ cannot be a set of lines through a
point, then $\phi_V(\Cal J_1)$ contains at least three non
collinear lines $\ell_1,\ell_2,\ell_3$. Then we can write
$\phi_V(\Cal J_1)=\ell_1 \cup \dots \cup \ell_s,$ where $s \geq 3$
and $\ell_i$ is a line ($1 \leq i \leq s$). Let $R$ be the
ramification divisor of $\phi_V$, then $R \in |K_S+3L|$. We can
write $R=R_1+\dots +R_s$ where $\phi_V(R_i)=\ell_i$. Moreover,
$R_i=\Sigma_{j=1}^{s_i}\alpha_{ij}R_{ij}$ where each $R_{ij}$ is
an irreducible curve and $\alpha_{ij} \geq 1$. Since $\phi_V$
ramifies along $R$ then for $1 \leq i \leq s$ there exists a
divisor $H_i=\phi_V^*(\ell_i) \in |V|$ and a divisor $E_i \geq 0$
such that $H_i=\Sigma_{j=1}^{s_j}(\alpha_{ij}+1)R_{ij}+E_i.$ We
claim that

\noindent (5.3.1) for any $1 \leq i \leq s$ there exists $j\ne i$
such that $R_i R_j>0$.

If (5.3.1) holds then for $x_{ij} \in R_i \cap R_j$ we get
$|V-R_i|,|V-R_j| \subset |V-2x_{ij}|$ and $|V-R_i|\ne |V-R_j|$.
Whence $x_{ij} \in \Cal J_2$ and codeg$(S,V) \geq 3$. Let us prove
(5.3.1). We argue with $R_1$ and the same argument holds when $i
\ne 1$. Suppose $R_1  R_j=0$ for any $j \ne 1$. Since $L$ is ample
then $0<L R_{11}=H_jR_{11}$ and so $$R_{11} E_{j}>0 \;{\text {for
any}}\; j \ne 1. \tag 5.3.2$$ This in particular implies $E_{i}>0$
for $1 \leq i \leq s$. Moreover, since $H_1$ is ample, its support
is connected and then
$$R_{11}(R_{12}+\dots +R_{1s_1}+E_1)>0. \tag 5.3.3$$ Now we have:
$$-(K_S+R_{11}) R_{11}=-(R_{11}+R-3L)R_{11}=-(R_{11}+R-H_1-H_2-H_3)R_{11}=$$
$$=-(R_{11}+R_1-H_1+R_2-H_2+R_3-H_3+\Sigma_{j \geq 4}R_j)R_{11}=$$
$$=(R_{12}+\dots +R_{1s_1}+E_1) R_{11}+
(R_{21}+\dots+R_{2s_2}+E_2) R_{11}+$$ $$
+(R_{31}+\dots+R_{3s_3}+E_3) R_{11}.$$ Note that the first summand
in the final expression is $\geq 1$ by (5.3.3) and the same
inequality holds for the other two summands by (5.3.2). Thus, by
adjunction formula we get:
$$-2 \leq 2g(R_{11})-2=(K_S+R_{11})R_{11} \leq-3,$$
a contradiction. This proves (5.3.1). Now suppose that
codeg$(S,V)=3$ (i.e, $s=3$ in the previous notation) and $E_1>0$.
Let $C_1$ be an irreducible component of $E_1$. Recall that $E_1$
has no non-reduced components (otherwise they would be part of
$R$). Then
$$-(K_S+C_1) C_1=- (R-3L+C_1)C_1=(R_{11}+\dots +R_{1s_1}+(E_1-C_1))
C_1+$$ $$+(R_{21}+\dots R_{2s_2}+E_2)
 C_1+(R_{31}+\dots R_{3s_3}+E_3) C_1\geq 3,$$
each of the three summands being $\geq 1$: the first one by the
connectedness of $H_1$ and the remaining two by the the ampleness
of $L$. By adjunction formula this gives a contradiction again.
This shows that $E_1=0$ and the same argument gives $E_2=E_3=0$.
Hence
$$K_S=-(R_{11}+\dots+R_{1s_1})-(R_{21}+\dots+R_{2s_2})-(R_{31}+
\dots+R_{3s_3}).$$ We claim that $s_1=s_2=s_3=1$. We have:
$-(K_S+R_{11}) R_{11}=(R_{12}+\dots+R_{1s_1})
R_{11}+(R_{21}+\dots+R_{2s_2}) R_{11}+(R_{31}+\dots+R_{3s_2})
R_{11}.$ If $s_1>1$ then by the connectedness of $H_1$ and the
ampleness of $L$ we have the contradiction $(K_S+R_{11})
R_{11}\leq -3$. The same argument works with $s_2$ and $s_3$.

Since $s_1=s_2=s_3=1$, $H_i=(\alpha_{i1}+1)R_{i1}$, $1 \leq i \leq
3$. Then $R_{i1}$ is an ample divisor for $1 \leq i \leq 3$.
Moreover $K_S=-R_{11}-R_{21}-R_{31}$ and so $S$ is a Del Pezzo
surface with $-K_S$ being the sum of three ample divisors. Hence
$S=\Bbb P^2$, $R_{11},R_{21},R_{31} \in |\Cal O_{\Bbb P^2}(1)|$
and $L=\Cal O_{\Bbb P^2}(\alpha_{11}+1)$. \qed
\enddemo

Let us observe that (5.3) gives a different proof of (3.5) in the
case of surfaces. Moreover, recall that for any $V \subseteq
H^0(S,L)$, $\Cal D(S,V)=\Cal D(S,H^0(S,L)) \cap |V|$ (possibly
set-theoretically). We know that $\text {codeg}(\Bbb P^2,H^0(\Bbb
P^2,\Cal O _{\Bbb P^2}(2)))=3$ and codeg$(\Bbb P^2,H^0(\Bbb
P^2,\Cal O _{\Bbb P^2}(r)))>3$ for $r \geq 3$. Hence, in the
previous discussion, either $\alpha_{11}=1$, $(S,L)=(\Bbb P^2,\Cal
O_{\Bbb P^2}(2))$ and $V \subset H^0(\Bbb P^2,\Cal O_{\Bbb
P^2}(2))$ is as in {\rm (d.4.2)} of $(3.4)$ or special projections
of the $r$-Veronese embedding of $\Bbb P^2$ have to be considered:
$V=\langle x_0^r,x_1^r,x_2^r \rangle$ provides an example of
codegree 3 for any $r$.

\proclaim{(5.4) Proposition} Let $(S,L,V)$ be as in $(0.0)$,
$\dim(S)=2$. If {\rm codeg}$(S,V)\leq 3$ and $\Cal J_2=\emptyset$
then $S$ is a ruled surface. \endproclaim

\demo{Proof} Choose a general vector subspace $V' \subseteq V$
such that $\dim(V')=3$. Then $\phi_{V'}(S)=\Bbb P^2$ and
codeg$(S,V') \leq 3$. Moreover $\cal J_2(V')=\emptyset$ because,
if not, $\Cal D(X,V')$ has a linear component of maximal dimension
and so $\Cal D(X,V)$ has a linear component of maximal dimension
contradicting $\Cal J_2(V)=\emptyset$. Hence, by the usual
expression of $c_2(J_1(L))$ in terms of the invariants of $S$, see
for example \cite{LPS1, A.1.1}, and (3.2.7) we have:
$$c_2(J_1(L))=e(S)+2K_S  L +3L^2 \leq \text{
codeg}(S,V')(L^2-1), \tag 5.4.1$$ where $e(S)$ is the topological
Euler--Poincar\'e characteristic of $S$. In particular $e(S)+2K_S
 L \leq -3.$ Hence either $e(S)<0$ or $2K_S  L<0$ and we
are done either by Castelnuovo--De Franchis theorem or by Enriques
theorem, see \cite{B}.  \qed
\enddemo

\head 6 Codegree 2 surfaces
\endhead

Now we deal with codegree two surfaces. Let us consider $(S,L,V)$
as in (0.0) with $\dim(S)=2$ and codeg$(S,V)=2$. By \cite{LPS1,
Thm. 2.8}, $\dim (\Cal D(S,V))=N-1$ and by (5.3), $\Cal D(S,V)$
has just one maximal dimensional irreducible component, say $\Cal
D^0$. Whence $\Cal J_2=\emptyset$ and $\Cal D^0$ is either an
irreducible quadric cone, or  a smooth quadric. We can exclude the
first possibility: If $\phi_V(S)^\vee$ is a cone then $\phi_V(S)$
would be degenerate and so $\phi_V(S)=\Bbb P^2$. Whence $\Cal D^0$
cannot be an irreducible quadric cone. To deal with the second
possibility note that if $\Cal D^0$ is a smooth quadric then, by
biduality, $\Cal D^{0\vee} \subset \phi_V(S)$. Hence, either
$\phi_V(S) \subset \Bbb P^3$ is a smooth quadric or
$\phi_V(S)=\Bbb P^2$ and there is a smooth conic in $\phi_V(\Cal
J_1)$. We will need the following general fact.

\proclaim{(6.1) Lemma} Let $(S,L,V)$ be as in $(0.0)$ with
$\dim(S)=2$ and $N \geq 3$. If there exists an
$(N-1)$--dimensional component $\Cal D^0 \subseteq \Cal D$ which
is a smooth quadric and any other irreducible components of $\Cal
D$ is linear then either $(S,L)$ is a scroll or $|\phi_V(\Cal
J_2)| \geq 2$.
\endproclaim

\demo{Proof}  By hypothesis $N \geq 3$ so $\phi_V(S) \subset \Bbb
P^3$ is a smooth quadric. Since the other components of $\Cal D$
are linear, $\phi_V(\Cal J_1)$ is a union of lines. Let us observe
that if $(S,L)$ is not a scroll then the ramification divisor $R
\in |K_S+2L|$ of $\phi_V$ is an ample and effective divisor, see
\cite{LP1, Thm. 2.5}. Consider $\ell_i \in \Cal T_i$ ($i=1,2$) a
general line in the ruling $\Cal T_i$ of $\phi_V(S)$. Then
$C_i=\phi_V^{-1}(\ell_i)$ is a smooth curve. Let us observe that
$C_1+C_2 \in |V|$, $C_1^2=C_2^2=0$ and $C_1 C_2=L^2/2$. Since
$(S,L)$ is not a scroll note that for any irreducible component
$C$ of $C_1$ or $C_2$ the branch locus of $\phi_{V|_C}:C \to \Bbb
P^1$ can be neither empty nor a point (see (2.3)). Then for any
component $C$ of $C_1$ or $C_2$ the restriction $\phi_V|_C:C \to
\Bbb P^1$ branches in at least two points. This means that $R$ has
at least four components, say $R_1$, $R_2$, $R_3$ and $R_4$, such
that $R_1  C_1>0$, $R_2 C_1>0$, $R_3  C_2>0$ and $R_4 C_2>0$.
Since any component of $R$ maps onto a line on $\phi_V(S)$ we have
$R_i^2 \leq 0$, $R_1  R_2=0$ and $R_3 R_4=0$. Moreover
$\phi_V(R_1)\ne \phi_V(R_2)$ and, since $R$ is ample, $R R_i>0$.
Hence there exist two components $R_1'$ and $R_2'$ of $R$ such
that $R_{1} R_{1}'>0$ and $R_{2} R_{2}'>0$. Now take $p \in R_{1}
\cap R_{1}'$ and respectively $q \in R_{2}\cap R_{2}'$. We claim
that $p,q \in \Cal J_2$ and this proves the lemma because
$\phi_V(p) \ne \phi_V(q)$. Since $p \in R_{1} \subset \Cal J_1$
then $|V-R_{1}| \subset |V-2p|$ and equivalently $|V-R_{1}'|
\subset |V-2p|$. Then $|V-2p|$ contains two different lines and so
$|V-2p|=|V-p|$, that is, $p \in \Cal J_2$. The same argument can
be applied to $q$. \qed \enddemo

Let us come back to the codegree two case. Suppose $\phi_V(S)
\subset \Bbb P^3$ to be a smooth quadric. By (6.1) $(S,L)$ is a
scroll over a smooth curve $B$ and $L\equiv C_0+bf$. We have
$L_{|f}=\Cal O_{\Bbb P^1}(1)$ for any fibre $f$. In particular
$\phi_V(f)$ is a line on $\phi_V(S)$ and all lines image by
$\phi_V$ of fibres of the scroll are on the same ruling $\Cal T_1$
of the quadric $\phi_V(S)$. Consider two general lines $\ell_1$
and $\ell_2$ on the other ruling $\Cal T_2$ of $\phi_V(S)$. It
holds that $\phi_V^{-1}(\ell_1)=C_1$ and
$\phi_V^{-1}(\ell_2)=C_2$, with $C_1$ and $C_2$ smooth irreducible
curves which are sections of the scroll. Moreover $C_1 \cap
C_2=\emptyset$. Thus, by \cite{Ha, Ex.\ V.2.2}, we have $S=\Bbb
P(\Cal E)$ where $\Cal E$ splits as $\Cal E=\Cal O_B\oplus \Cal L$
with $\Cal L \in$ Pic $(C)$ and deg$(\Cal L) \leq 0$. Consider
$\ell$ a general line in the ruling $\Cal T_1$. Since
$\phi_V^{-1}(\ell_1 \cup \ell)$ is linearly equivalent to
$\phi_V^{-1}(\ell_2 \cup \ell)$ then $C_1$ and $C_2$ are linearly
equivalent. In particular this says that $S=B \times \Bbb P^1$ and
we exactly are in case (b.1) of (1.1.1) with $n=2$.

\noindent By what is said just before (6.1) it remains to consider
the case when $\phi_V(S)= \Bbb P^2$ and $\phi_V(\Cal J_1)$
contains a smooth conic $C$. In this situation $\Cal D(S,V)=C^\vee
\cup \Cal D^1 \cup \dots \cup \Cal D^s$ where, for $1 \leq i \leq
s$, each $\Cal D^i$ is a point in $|V|$ corresponding to a line in
$\phi_V(\Cal J_1)$.

Suppose for the moment that $(S,L)$ is a scroll. Take a general
$\ell_1 \in C^\vee$ corresponding to $D_1 \in \Cal D(S,V)$. By
(2.5.4) for any $x \in \text{Sing}(D_1)$ it holds that
$D_1=f_{\pi(x)}+R$ with $R$ an effective divisor smooth at $x$. In
particular $\mu_x(D_1)=1$ and $c_2(J_1(L))=L^2$ is an even number
by (3.2.3). Then, by the same argument as in (2.5.4), we obtain
$D_1=M+f_1+\dots+f_{L^2/2}$. Since $L^2=D_1^2$ then $M^2=0$. The
same construction can be done for another general line $\ell_2 (
\ne \ell_1) \in C^\vee$. In particular
$D_2=M'+f'_1+\dots+f_{L^2/2}'$. Now $\phi_V(M)=\ell_1 \ne
\ell_2=\phi_V(M')$ and so $M \ne M'$. Since $M  M'=0$ and $M$,
$M'$ are irreducible, we have constructed a one dimensional family
of pairwise disjoint sections. Moving the point $\ell$ on $C^\vee$
we get a rational parametrization of these sections $M$. So $M$ is
linearly equivalent to $M'$ showing that we again obtain a product
of a smooth curve cross $\Bbb P^1$ as in (b.1) of (1.1.1) with
$n=2$.

Suppose now that $(S,L)$ is not a scroll. In this case (3.2.7)
reads as $e(S)+2K_SL+3L^2 \leq 2L^2-2$, where $e(S)$ is the
topological Euler--Poincar\'e characteristic of $S$. In particular
this gives:
$$e(S)+2K_SL \leq -L^2-2. \tag 6.2$$ Since $(S,L)$ is neither a
scroll nor $(\Bbb P^2, \Cal O_{\Bbb P^2}(r))$ with $r=1,2$ then
$(K_S+L)^2 \geq 0$ by \cite{LP1, 2.1}, or equivalently:
$$2K_SL \geq -K_S^2-L^2. \tag 6.3$$ Substituting (6.3) in (6.2) we have:
$$-L^2-2 \geq e(S)+2K_SL \geq e(S)-K_S^2-L^2. \tag 6.4$$ In
particular $$e(S)-K_S^2 \leq -2. \tag 6.5$$ By (2.5.4) $S$ is
ruled, then either $(S,L)=(\Bbb P^2, \Cal O_{\Bbb P^2}(r))$, $r
\geq 3$, or $S \ne \Bbb P^2$, $e(S)=4(1-q)+s$ ($s$ is a
nonnegative integer), $K_S^2=8(1-q)-s$. By Bezout's theorem

\noindent (6.6)  the sum of the Milnor numbers of the
singularities of a plane curve of degree $r$ with isolated
singularities is less than or equal to $(r-1)^2$.

If $(S,L)=(\Bbb P^2, \Cal O_{\Bbb P^2}(r))$ then by (6.6), (3.2.7)
and (3.3.2) we get the contradiction $3(r-1)^2 =c_2(J_1(L)) \leq
2(r-1)^2$.

Let us observe the following general fact: take a general $D \in
\Cal D$ which corresponds to a line tangent to $C$ at $y$. Then,
using (3.2) and (3.2.6) we get:

\noindent (6.7) If $c_2(J_1(L))=2(L^2-1)$ then
$\phi_V^{-1}(y)=\{x\} \subset$Sing$(D)$ and $\mu_x(D)=L^2-1$.

\noindent (6.8) If $c_2(J_1(L))=2(L^2-2)$ then
$\phi_V^{-1}(y)=\{x,x'\} \subset$Sing$(D)$ and
$\mu_x(D)+\mu_{x'}(D)= L^2-2$ (possibly $x=x'$ and $\mu_{x}(D)=
L^2-2$).

\noindent If $(S,L) \ne (\Bbb P^2, \Cal O_{\Bbb P^2}(r))$ then
(6.4) gives
$$-L^2-2 \geq 4(1-q)+s+2K_SL\geq
4(1-q)+s-K_S^2-L^2\geq-4(1-q)+s-L^2.$$ In particular $s+2 \leq
4(1-q)$, that is, $q=0$ and $s \leq 2$. Hence, by Noether formula,
we have $K_S^2=8-s$. Then either $q=0$, $s=2$, $K_S^2=6$,
$e(S)=6$, not compatible with (6.5); or $q=0$, $s=1$, $K_S^2=7$,
$e(S)=5$; or $q=0$, $s=0$, $K_S^2=8$, $e(S)=4$.

\noindent If $q=0$, $s=1$, $K_S^2=7$, $e(S)=5$ then equality holds
in (6.4), hence in both (6.2) and (6.3). By \cite{LP1,2.1}
equality $(K_S+L)^2=0$ implies that either $S$ is a Del Pezzo
surface, $L \equiv -K_S$, $c_2(J_1(L))=12$, (excluded by (6.6)
because we are dealing with plane cubic curves), or

\noindent (6.9) $(S,L)$ is a rational conic bundle, more
precisely, $S$ is a blowing up at $s=1$ point of a $\Bbb P^1$
bundle over $\Bbb P^1$ of invariant $e \geq 0$. Denote by $E$ the
exceptional divisor and $C_0$ and $f$, as an abuse of notation,
the proper transforms of the corresponding $C_0$ and $f$ on the
$\Bbb P^1$ bundle. Hence $L=2C_0+bf-E$, $L|_f=\Cal O_{\Bbb
P^1}(2)$, $c_2(J_1(L))=2(L^2-1)$ and $L^2=4(b-e)-1$.

\noindent If $q=0$, $s=0$, $K_S^2=8$, $e(S)=4$ then $S$ is a
rational $\Bbb P^1$--bundle of invariant $e \geq 0$ and $(S,L)$ is
not a scroll. Whence $L=aC_0+bf$, $a \geq 2$ and $b>ae$. By
(5.4.1) we have
$$4+2ae-4a-4b+6ab-3a^2e \leq 4ab-2a^2e-2, \tag 6.9.1$$ or equivalently
$b(2a-4)-a^2e+2ae-4a \leq -6.$ Since $a \geq 2$ and $b \geq ae+1$
then $(ae+1)(2a-4)-a^2e+2ae-4a \leq -6.$ Then $a(e(a-2)-2) \leq
-2$ which gives that either $e=0$, or $e=1$, $a=2,3$ or $e \geq 2$
and $a=2$. If $a=2$ then we get (6.10) else we get the following.
In the case $e=0$ we can suppose $0<a \leq b$ and the inequality
(6.9.1) implies that $a=3=b$. If $e=1$ and $a=3$ then
$c_2(J_1(L))$ is odd, contradicting (3.2.3). If $e=0$, $a=b=3$
then the ramification divisor's class is $R=K_S+3L=7C_0+7f$ and
$c_2(J_1(L))=2(L^2-1)=34$. By (6.7), moving the singular point,
there exists an effective divisor $F$ on $S$ such that $R-17F \geq
0$. This clearly gives a contradiction.

\noindent (6.10) If $a=2$ then $S$ is a $\Bbb P^1$ bundle over
$\Bbb P^1$, $L|_f=\Cal O_{\Bbb P^1}(2)$ for any fibre $f$ of $S$,
$c_2(J_1(L))=2(L^2-2)$ and $L=2C_0+bf$. We just need to face (6.9)
and (6.10).

Let us observe in (6.10) that ampleness is equivalent to very
ampleness, hence $\phi_V(S)$ is just the projection of $\phi_L(S)
\subset \Bbb P^{h^0(S,L)-1}$ from a codimension three linear space
$T$ such that $T \cap \phi_L(S)=\emptyset$. A similar situation
occurs in (6.9). We can blow down the exceptional divisor to
obtain $S'$. Consider the line bundle $L'=2C_0+bf$ that is, in
fact, very ample. Take $V' \subset H^0(S',L')$ defining a linear
system with just one base point.  Then $\phi_V$ is the morphism
resolving the indeterminacy of the rational map defined by $|V'|$.
Hence we are projecting $\phi_{L'}(S') \subset \Bbb
P^{h^0(S',L')-1} $ from a codimension 3 linear space $T$ meeting
$\phi_{L'}(S')$ in one point.

First let us deal with the case (6.10). In this case $L^2=4(b-e)$
and $c_2(J_1(L))=2(L^2-2)>L^2$ since $(S,L)$ is not a scroll,
hence
$$L^2=4(b-e)\geq 8. \tag 6.10.1$$ Consider a general $D \in |V|$.
By Riemann--Hurwitz formula the ramification divisor $R_D$ of
$\phi_{V|D}$ verifies deg$(R_D)=L(L+K_S)+2L^2=10(b-e)-4.$ On the
other hand there exist effective divisors $F_i=a_i C_0+b_i f >0$
($1 \leq i\leq s$) such that $|V-2F_i| =\Cal D^{i}$ and
$\phi_V(F_i)=\ell_i$. Then there exist integers $\alpha_i \geq 1$
and a divisor $G>0$ such that the ramification divisor $R$ of
$\phi_V$ verifies $R=\alpha_1 F_1+\dots +\alpha_s F_s +G$ and
$\phi_V(G)=C$. Moreover, by (6.8), there exist two divisors
$G_1=A_1C_0+B_1f>0, G_2=A_2C_0+B_2f>0$ (maybe equal) and two
integers $z_1,z_2 \geq 0$, $z_1+z_2=L^2-2$ such that
$$R=\alpha_1 F_1+\dots +\alpha_s F_s+z_1G_1+z_2G_2.\tag 6.10.2$$
Since any ramification
point of $\phi_{V|D}$ is a ramification point of $\phi_V$ we get
$$10(b-e)-4=\deg(R_D)=RL \leq \alpha_1+\cdots +\alpha_s+L^2-2=
\alpha_1+\dots +\alpha_s+4(b-e)-2,$$ which implies
$$\alpha_1+\dots +\alpha_s \geq 6(b-e)-2. \tag 6.10.3$$ By the
ramification formula $R=K_S+3L=4C_0+(3b-2-e)f$ then:
$$\matrix
\alpha_1a_1+\dots+\alpha_sa_s+z_1A_1+z_2A_2=4\phantom{b-e-2.} \cr
\alpha_1b_1+\dots+\alpha_sb_s+z_1B_1+z_2B_2=3b-e-2. \endmatrix$$
Note that the effectiveness of $F_i$ and $G_i$ implies that
$a_i+b_i \geq 1$ and $A_i+B_i \geq 1$. Adding both equalities and
using (6.10.3) we get $7b \leq 9e+6$. This is a contradiction
because $b\geq 2e+1$ by the ampleness of $L$.

If (6.9) holds we can argue in exactly the same way. In fact
$L^2=4(b-e)-1$ and $c_2(J_1(L))=2(L^2-1)$ then the expression for
the ramification divisor is $$R = \alpha_1 F_1 + \dots + \alpha_s
F_s + z_1G=K_S+3L=4C_0+(3b-2-e)f-2E \tag 6.10.4$$ with
$z_1=L^2-1$, $F_i=a_iC_0+b_if+c_iE$ and $G=AC_0+Bf+CE$. The
formula now gives $\deg(R_D)=10(b-e)-6$. Whence the analogue of
(6.10.3) is $\alpha_1+\dots+\alpha_s\geq 6(b-e)-4$. Since $G>0$
and $F_i>0$ ($1\leq i \leq s$), $A+B+C \geq 1$ and $a_i+b_i+c_i
\geq 1$. Then, by (6.10.4) and the previous inequality we get $7b
\leq 9e+6$, a contradiction.

Summing up the discussion in the codegree two case we have proved
the following

\proclaim{(6.11) Theorem} Let $(X,L,V)$ be as in $(0.0)$ with
$\dim(X)=2$ and \text{\rm codeg}$(X,V)=2$. Then $X=C \times \Bbb
P^1$ for a smooth curve $C$, $|V_{|C}|$ is a $g_d^1$ on $C$
defining a degree $d$ morphism $f:C \to \Bbb P^1$, and $L=F^*\cal
O_{Q}(1)$ where $F=s \circ (f \times Id)$ is the composition of $f
\times Id: C \times \Bbb P^1 \to \Bbb P^1 \times \Bbb P^1$ with
the Segre embedding $s:\Bbb P^1 \times \Bbb P^1 \to Q \subset \Bbb
P^3$.
\endproclaim

\head 7 The codegree 3 case.
\endhead

Let $(S,L,V)$ be a triplet as in (0.0) such that $\dim(S)=2$ and
codeg$(S,V)=3$. Let $\Cal D^0$ be the union of the
$(N-1)$--dimensional irreducible components of $\Cal D(S,V)$. We
have already considered the case of $\Cal D^0$ being the union of
three distinct hyperplanes in (5.3). Hence one of the following
holds:

(7.1) $\Cal D^0=Q \cup H$, where $Q$ is a smooth quadric
hypersurface and $H=\Bbb P^{N-1}$,

(7.2) $\Cal D^0=\phi_V(S)^\vee$,

(7.3) $\Cal D^0=Q \cup H$, where $Q$ is an irreducible quadric
cone and $H=\Bbb P^{N-1}$,

(7.4) $\Cal D^0=\phi_V(C)^\vee$ for $C \subset \Cal J_1$ an
irreducible curve.

Consider (7.1) Since $Q^\vee \subseteq \phi_V(S)$ then either
$\phi_V(S) \subset \Bbb P^3$ is a smooth quadric hypersurface or
$\phi_V(S)=\Bbb P^2$ and $Q$ is a smooth conic. If $\phi_V(S)$ is
a smooth quadric hypersurface then, by (6.1), either $(S,L)$ is a
scroll, so that $\Cal J_2=\emptyset$, or $|\phi_V(\Cal J_2)| \geq
2$. Hence this case does not occur. It remains to consider
$\phi_V(S)=\Bbb P^2$ and $\phi_V(\Cal
J_1)=C\cup\ell_1\cup\dots\cup \ell_s$, where $C$ is a smooth conic
and $\ell_i$ is a line for $1 \leq i\leq s$. This case is
effective as shown in (d.4.1) of (3.4).

Now consider (7.2). Since $\dim(\phi_V(S)^\vee)=N-1$, $\phi_V(S)$
cannot be a cone. By (5.1), either $(S,L)=(\Bbb P^2, \Cal O_{\Bbb
P^2}(2))$ (discussed in (d) of (3.4)) or deg$(\phi_V(S)) \leq 3$
and $(S,L)$ is a scroll.

When $(S,L)$ is a scroll and deg$(\phi_V(S))=3$  then $N=3,4$.
Recall that $\phi_V(S)$ is not a cone. If $N=4$ then $\phi_V(S)
\subset \Bbb P^4$ is the rational normal scroll $\Bbb P(\Cal
O_{\Bbb P^1}(-1) \oplus \Cal O_{\Bbb P^1}) \subset \Bbb P^4$,
\cite{XXX}, which has a (single) section that is a line, say $C$.
If $N=3$ then $\phi_V(S)$ is a (finite and birational) projection
of $\Bbb P(\Cal O_{\Bbb P^1}(-1) \oplus \Cal O_{\Bbb P^1}) \subset
\Bbb P^4$ from a point $p \in \Bbb P^4 \setminus \Bbb P(\Cal
O_{\Bbb P^1}(-1) \oplus \Cal O_{\Bbb P^1})$. In view of the
classification of cubic ruled surfaces \cite{E, Ch. I,\ \S\S 37,
38, 48, 49} there are two types:

(7.2.1) $p$ lies on the plane spanned by a smooth conic $Q \subset
\Bbb P(\Cal O_{\Bbb P^1}(-1) \oplus \Cal O_{\Bbb P^1}) \subset
\Bbb P^4$. Then the image by the projection from $p$ has two
directrix lines $C_1$ and $C_2$ (lines meeting each line of the
ruling), being $C_1$ the projection of $C$ and $C_2$ the
projection of $Q$ (a double line).

(7.2.2) $p$ lies on the plane plane spanned by a general fiber $f$
and $C$. Then the image by the projection has just one directrix
line $C_1$ that is the projection of $C$ (or of $f$, so that a
double line).

Let us consider first $N=4$. Take a general point $p \in
\phi_V(S)$ and a general hyperplane section $H$ singular at $p$.
Then $H$ decomposes as the fibre through $p$, say $f_p$ and a
conic $Q_p$ cutting $f_p$ just at $p$. Since
$c_2(J_1(L))=L^2=$codeg$(X,V)$ and $\deg(\phi_V(S))=3$ then, by
(3.2.4), $\phi_V^*(H)$ is singular at exactly the $L^2/3$ points
constituting $\phi_V^{-1}(p)$. Then, by (2.5.4),
$\phi_V^{*}(H)=M+f_1+\dots +f_{L^2/3}$ where $M$ is a smooth curve
and $f_i$ is a fibre of $S$ such that $\phi_V(f_i)=f_p$. In
particular $\phi_V^*(f_p)=f_1+\dots +f_{L^2/3}$. Take now $H$ cut
out by the hyperplane containing two general fibres $f_q$ and
$f_p$. Then $H=C+f_p+f_q$. By the previous arguments
$\phi_V^*(H)=D+f_1+\dots+f_{L^2/3}+g_1+\dots+g_{L^2/3}$ where
$D=\phi_V^{*}(C)$ is a curve such that $D^2=-L^2/3$ and
$\phi_V^*(f_q)=g_1+\dots +g_{L^2/3}$. Now $D
M=D(D+g_1+\dots+g_{L^2/3})=0$. This implies that $\Cal E$ is
decomposable \cite{Ha, Exercise 2.2 p.\ 383}. Assume that $\Cal E$
is normalized, in the usual sense \cite{Ha, p.\ 373}; then $\Cal E
= \Cal O_B \oplus \Cal L$, where $\Cal L \in \text{Pic}(B)$ is
such that $e:=-\deg \Cal L \geq 0$ \cite{Ha, Theorem 2.12 p.\
376}. Let $C_0$ be, as usual, a tautological section. We can write
$D \equiv C_0 + af$ for some integer $a$ and $M \equiv C_0 + bf$,
where $b=a+\frac{L^2}{3}$. Moreover, $0 = D  M= -e + a + b$. As
$D$ is irreducible, we know from \cite{Ha, Proposition 2.20 p.\
382} that either $D=C_0$ or $a \geq e$. However, in the latter
case we get $e = a+b = 2a + \frac{L^2}{3} \geq 2e +
\frac{L^2}{3}$, giving $e + \frac{L^2}{3} \leq 0$, a
contradiction. Therefore $D=C_0$, hence $a=0$ and then
$\frac{L^2}{3} = b = e$. In particular, $\phi_{V|D}:D \to B$ is an
isomorphism and $L \equiv C_0 + 2 \frac{L^2}{3} f$.

If $N=3$, as said before, there exists a finite and birational
morphism $\pi_p: \Bbb P(\Cal O_{\Bbb P^1}(-1) \oplus \Cal O_{\Bbb
P^1}) \to \phi_V(S)$ (the projection from $p$) which is the
normalization morphism. Then, by the universal property of the
normalization, $\phi_V:S \to \phi_V(S)$ factors through $\pi_p$.
Now we can get the same conclusion as in the previous paragraph
with the following warning. For the general line $f_p$ on
$\phi_V(S)$, we can consider
$\phi_V^*(f_p)=\pi^{*}(f_p)=f_1+\dots+f_{L^2/3}$ and $M$ is
defined exactly as before. To construct $D$ we can do the
following. In (7.2.1) just consider a hyperplane section $H$
containing $C_1$ and $f_p$, then there exists another line of the
ruling, say $f_q$ and $H=C_1+f_p+f_q$ and proceed as before. In
(7.2.1) take a hyperplane section containing $C_1$ and $f_p$. Then
$\phi_V^*(H)=\pi^*(C+f+f_p)$ and proceed as before.

If deg$(\phi_V(S))=2$ then either $\phi_V(S)$ is a quadric cone or
$\phi_V(S) \subset \Bbb P^3$ is a smooth quadric. If the former
occurs then we get the contradiction that $\Cal D^0$ is
degenerate. The latter contradicts deg$(\phi_V(S)^\vee)=3$.

Next consider (7.3). If $\Cal D^0$ contains an irreducible quadric
cone $Q$ then by (1.4) $Q^\vee \subset \phi_V(S)$. Hence $Q^\vee$
is a smooth plane conic. Moreover dim$(\phi_V(S))^\vee <N-1$, and
so, by (5.2), $N=3$ and $\phi_V(S)$ is a quadric cone.

Finally we deal with (7.4). Since $\dim(\phi_V(S)^\vee)<N-1$ and
$\Cal J_1=\emptyset$ then, by (5.2), $\phi_V(S)=\Bbb P^2$.

Summing up the discussion on the codegree three case we get:

\proclaim{(7.5) Theorem} Let $(X,L,V)$ as in $(0.0)$ with
$\dim(X)= 2$ and \text{\rm codeg}$(X,V)=3$. Then, either

$(7.5.1)$ $(X,L)=(\Bbb P^2, \Cal O_{\Bbb P^2}(r))$, or

$(7.5.2)$ $X=\Bbb P_B(\Cal O_B \oplus \Cal O_B(p^*\Cal O_{\Bbb
P^1}(-1))) \mapright{\pi} \; B$, where $B$ is a smooth curve, $p:B
\to \Bbb P^1$ is a surjective morphism and $L\equiv
C_0+\pi^*p^*(\Cal O_{\Bbb P^1}(2))$; $N=3,4$ and $\phi_V(X)$ is a
cubic ruled surface, or

$(7.5.3)$ $N=3$, $\phi_V(X) \subset \Bbb P^3$ is an irreducible
quadric cone, $\phi_V(\Cal J_1)=C \cup \ell_1 \cup \dots \cup
\ell_s$ where $C$ is a smooth plane conic and $\ell_i$ is a line
for $1 \leq i \leq s$, or

$(7.5.4)$ $\phi_V(X)=\Bbb P^2$ and there exists an irreducible
curve $C \subseteq \Cal J_1$ such that $\phi_V(C)^\vee$ is one of
the maximal dimensional components of $\Cal D$.
\endproclaim

In (3.4.d) we have studied (7.5.1) when $r=2$; (3.3.d) with $n=2$
also provides examples of (7.5.1) for any $r \geq 3$. Let us
observe that the general hyperplane section of the Segre variety
$\Bbb P^1 \times \Bbb P^2 \subset \Bbb P^5$ is the only example as
in (7.5.2) in the classical setting, \cite{Z2, p. 93}. We end this
section with examples corresponding to (7.5.2), (7.5.3) and
(7.5.4). Then all situations described in (7.2) are effective.

\noindent{\bf (7.6) Examples}

{\bf (a)} Let $B$ be an elliptic curve, and let $p^{\prime}:B \to
\Bbb P^1$ be the morphism defined by a line bundle of degree $2$
on $B$. Let $p: S \to \Bbb F_1$ be the double cover branched along
the fibres of $\Bbb F_1$ corresponding to the four branch points
of $p^{\prime}$. Then $S$ is a $\Bbb P^1$ bundle over $B$.
Moreover its invariant is $2$. To see this denote by $\gamma_0$
the $(-1)$-section of $\Bbb F_1$ and note that
$C_0=p^{-1}(\gamma_0)=p^*\gamma_0$ is the section of minimal
self-intersection on $S$. Set $L:=p^*[\gamma_0+2\varphi]$, where
$\varphi$ is a fibre of $\Bbb F_1$, and $V:=p^*H^0(\Bbb F_1,
[\gamma_0+2\varphi])$. Let $f$ be the general fibre of $S$. As
$p^* \varphi$ consists of two fibres of $S$, we have $L \equiv C_0
+ 4f$. Note that $L$ is ample by \cite{Ha, Proposition 2.20 p.\
382}; moreover $V$ spans $L$ by construction, $L^2=6$ and
$\phi_V|_{C_0}=p$. This gives an example of (7.5.2) with $g(B)>0$.

{\bf (b)} Let $\nu: \Bbb F_2 \to \Gamma \subset \Bbb P^3$ be the
minimal desingularization of the quadric cone. Let $C_0$ and $f$
be the minimal section and a fibre of $\Bbb F_2$. Note that $\nu^*
(\Cal O_{\Gamma}(1)) = [C_0+2f]$. Let $C \in |C_0+2f|$ be a smooth
curve (the pull-back of a general hyperplane section $\gamma$ of
$\Gamma$). Then $\Delta:= C_0+C$ is a smooth divisor in the linear
system $|2\Cal B|$, where $\Cal B = [C_0+f]$. Thus there exists a
smooth surface $X$ and a morphism $\rho:X \to \Bbb F_2$ of degree
$2$ branched along $\Delta$. Let $E:=\rho^{-1}(C_0)$; thus
$\rho^*(C_0) = 2E$, since $C_0$ is in $\Delta$. Moreover, $4E^2 =
(2E)^2 = (\rho^*C_0)^2 = 2 C_0^2 = -4$. So $E$ is a $(-1)$-curve
inside $X$, and we can contract it, obtaining a smooth surface
$S$. Let $\mu:X \to S$ be the contraction and set $x=\mu(E)$. Then
we get a commutative diagram
$$\CD
  X   @>\mu>>  S  \\
  @V \rho VV    @V \pi VV   \\
  \Bbb F_2  @>\nu>>   \Gamma\ ,
  \endCD $$
where $\pi:S \to \Gamma$ is the induced double cover. Note that
$\pi$ is branched along $\gamma = \nu(C)$ and at  vertex $v$ of
$\Gamma$, and $v=\pi(x)$. Put $L:=\pi^*\Cal O_{\Gamma}(1)$ and $V
= H^0(S,L)$. Then $\phi_V = \pi$. We have $\Cal J_2(S,V)=\{x\}$,
while $\Cal J_1(S,V) \setminus \Cal J_2(S,V) = \pi^{-1}(\gamma)$.
It follows that $\Cal D(S,V)$ consists of: $\Cal D_0$, the dual of
$\Gamma$, which is a smooth conic; $\Cal D_1$, the dual of
$\gamma$, which is a quadric cone (because $\gamma$ is a plane
curve) containing $\Cal D_0$; and $\Cal D_2$, the plane of $\Bbb
P^{3 \vee}$ parameterizing the planes through the vertex $v$.
Therefore $\text{codeg}(S,V)= \deg \Cal D_1 + \deg \Cal D_2 = 3$.
This gives an example as in (7.5.3).

It deserves to explore the example above a little bit more, to
recognize a situation early described in (3.4). Note that the
ruling projection $\Bbb F_2 \to \Bbb P^1$ induces a fibration $X
\to \Bbb P^1$, whose general fibre $F:=\rho^*(f)$ is a $\Bbb P^1$,
being a double cover of $f$ branched at $\Delta \cap f$. Hence $X$
is rational, and so is $S$. By the ramification formula we have
$$K_X = \rho^*(K_{\Bbb F_2}+ \Cal B) = \rho^*(-2C_0-4f+C_0+f)=
- \rho^*(C_0+3f). \tag 7.6.1 $$ It thus follows that
$K_X^2=2(C_0+3f)^2=8$, and so $K_S^2 = K_X^2+1=9$. Therefore $S =
\Bbb P^2$. From the commutativity of the diagram above, recalling
(7.6.1) and the fact that $K_X = \mu^* K_S + E$, we also see that
$$\split
\mu^* (\pi^* \Cal O_{\Gamma}(1)) &= \rho^* (\nu^* \Cal
O_{\Gamma}(1)) = \rho^* (C_0+2f) = (2E + 2F) \\
&= \frac{2}{3}(2E+3F + E)=\frac{2}{3}(-K_X + E) = \frac{2}{3}\mu^*
(-K_S). \endsplit $$ Therefore $\pi^* \Cal O_{\Gamma}(1) =
\frac{2}{3}(-K_S) = \Cal O_{\Bbb P^2}(2)$. This shows that $(S,L)
= (\Bbb P^2, \Cal O_{\Bbb P^2}(2))$, and we fall in case (d.2) of
(3.4).

{\bf (c)} In (7.5.4), as $N=2$, either $\phi_V(C)$ is a smooth
plane conic and $|\phi_V(\Cal J_2)|=1$ or $\phi_V(C)^\vee=\Cal D$,
$\cal J_2=\emptyset$ and so $S$ is a ruled surface by (5.4). Note
that for a general codimension one subvector space $V' \subset V$
of any example as in (7.5.3) we get the former situation. On the
other hand something else can be said thanks to Pl\"ucker
formulas. If $\phi_V(C)$ has ordinary singularities then either

$\phi_V(C)$ is a sextic with nine cusps (and no other
singularities) and $\phi_V(C)^\vee$ is a smooth plane cubic (this
example is effective as shown in (d) of (2.6), and in (d.3.1) of
(3.4)), or

$\phi_V(C)$ is a quartic with three cusps (and no other
singularities) and $\phi_V(C)^\vee$ is a nodal cubic (this example
is effective as shown in (d.3.2) of (3.4)), or

$\phi_V(C)$ and so $\phi_V(C)^\vee$ are cuspidal cubics. Let us
put an example of this last situation. Take $X=\Bbb P^1 \times
\Bbb P^2$ and $L$ the line bundle defining the Segre embedding in
$\Bbb P^5$. Consider $V \subset H^0(X,L)$ defining a general
base-point free linear system with $\dim(|V|)=3$. Since
$X^\vee=\Bbb P^1 \times \Bbb P^2 \subset \Bbb P^{5 \vee}$ we have
that $X^\vee \cap |V|=C_0 \subset \Bbb P^3$ is a twisted cubic.
Call $TC_0 \subset \Bbb P^3$ the tangent variety to $C_0 \subset
\Bbb P^3$, that is, the union of its tangent lines. Consider $H
\in TC_0 \setminus C_0$. In particular $H \notin X^\vee$. Then $X
\cap H=S \subset \Bbb P^4$ is a smooth cubic scroll. Let us
observe that the restriction of $|V|$ to $S$ is a base-point free
linear system of dimension 2. In fact one can suppose that
$V=\langle s_0,s_1,s_2,s_3\rangle$, where $s_0$ defines $H$. Then,
when restricting to $S$, $V|_S=\langle s_1|_S, s_2|_S, s_3|_S
\rangle$. If there exists $x \in \text{Bs}|V|_S|$ then $s_0,
\dots, s_3$ vanish at $x$, contradicting that $|V|$ is base-point
free. It is classically well known that the projection $\pi_H$
from $H$ gives the identification $S^\vee =\Cal D(S,
L|_S)=\pi_H(X^\vee)=\pi_H(\Cal D(X,L))$; then $\Cal
D(S,V|_S)=\pi_H(\Cal D(X,L)) \cap |V|_S|=\pi_H(\Cal D(X,L)) \cap
|V|)=\pi_H(C_0)$ a cuspidal curve.
\medskip

\head 8 Final remarks
\endhead

In this section we present some problems which we consider of
interest.

\noindent {\bf (8.1)} As pointed out after (5.1), in the classical
setting it is possible to classify surfaces for which the
difference
$$c_2(J_1(L))-L^2=\text{\rm class}(S)-\text{deg}(S) \tag 8.1.1$$ is
small.  In the ample and spanned case, triplets $(S,L,V)$ for
which the right hand term of (8.1.1) is less than or equal to zero
are listed in \cite{LPS1, Prop.\ A.1}. In line with this we have
stated (5.1) where surfaces for which codeg$(S,V)-\deg(\phi_V(S))
\leq 0$ are considered. In this context it has sense the following
definition:

\proclaim{(8.1.2) Definition} Let $(X,L,V)$ be a triplet as in
$(0.0)$. We say that $(X,V)$ has {\rm tame codegree} if {\rm
codeg}$(X,V) =c_n(J_1(L))$.
\endproclaim

Pairs in examples (a) and (c) of (3.3) have tame codegree while
for (b) (and $d \geq 2$) in (3.3) we have codeg$(X,V)
<c_n(J_1(L))$. In the classical setting, i.\ e., when $\phi_V$
gives an embedding, having tame codegree simply means that $\Cal
D(X,V)$ is a hypersurface, because in that case
$c_n(J_1(L))=\deg(\Cal D(X,V))$. More generally, in the ample and
spanned setting, haveing tame codegree means that the general
element in $\Cal D$ is singular in a single point and the
singularity is just a non-degenerate quadratic singularity, see
(3.2.4). Let us show another example.

\flushpar $(8.1.3)$ Let $S=\Bbb P(\Cal E)\mapright{\pi} B$, where
$\Cal E$ is the rank-2 vector bundle over a smooth curve $B$ of
genus 1 of (2.6.d), defined by a non-split exact sequence
$$0 \to \Cal O_B \to \Cal E \to \Cal O_B(p) \to 0 \qquad (p \in
B),$$ and $L=2\xi$, where $\xi$ is the tautological line bundle of
$\Cal E$. Clearly $L$ is ample since $\xi$ is so. Moreover, $L$ is
spanned as Reider's Theorem immediately shows. Note that
$h^0(L)=h^0(S^2\Cal E)=\deg (\Cal E) =3$. So $(S,L,V=H^0(S,L))$ is
as in (0.0). Let us show that $(S,V)$ has tame codegree. We can
regard $S$ as the twofold symmetric product of the base elliptic
curve $B$. Hence $L$ lifts to $B \times B$ via the natural double
cover $p:B \times B \to S$ as the line bundle $p_1^* \Cal O_B(x+y)
\otimes p_2^* \Cal O_B(x+y)$, where $p_i: B \times B \to B$ is the
projection onto the $i$-th factor, and $x,y \in B$. On the other
hand, since $B$ is an elliptic curve, for any two points $x, y \in
B$ the linear series $|x+y|$ is a $g^1_2$. So our $L$ is like that
appearing in \cite{BDL, Ex. 9}, where the ramification locus of
$\phi_V$ is described. In fact the branch locus of the 4--to--1
map $\phi_V:S \to \Bbb P^2$ is the union of a smooth conic and
four of its tangent lines. For any $b \in B$ we will denote
$f_b=\pi^*\Cal O_B(b)$ and $\Cal O_S(C_0)=\xi$. Note that $L^2=4$,
hence $c_2(J_1(L))=8$. Observe that $\phi_V$ embeds $f_p$ as a
smooth plane conic, say $\gamma$, and gives a 2--to--1 map from
$C_0$ onto a line. For any $t \in B$ note that $h^0(S,\Cal
O_S(C_0+f_p-f_t))=h^0(S,\Cal O_S(C_0+f_t-f_p))=1$. Then we can
choose $\Gamma_t \in |C_0+f_p-f_t|$ and $\Gamma'_{t}\in
|C_0+f_t-f_p|$ such that $\Gamma_t+\Gamma'_{t} \in |2C_0|=|L|$ and
so $\phi_V(\Gamma_t)=\phi_V(\Gamma'_{t})$. Moreover one can check
by \cite{BDL, Ex. 9} that $\phi_V(\Gamma_t+\Gamma'_{t})$  meets
$\gamma$ in just one point. This gives $\gamma^\vee \subseteq \Cal
D(S,V)$. Moreover there exist $p_i \in B$, $i=1,2,3$, such that
$\Cal O_B(2p_i)=\Cal O_B(2p)$. Call $p_0=p$. This produces four
non--reduced elements $2\Gamma_{p_i} \in |2C_0|$. Hence the six
lines $\langle 2\Gamma_{p_i},2\Gamma_{p_j}\rangle \subset |V|$ are
contained in $\Cal D(S,V)$. This gives codeg$(X,V)\geq 8$ and in
fact an equality by (3.2.2).

In the following paragraphs we classify surfaces as in (0.0)
having tame codegree $\leq 8$. The argument relies on some rough
inequalities. In fact, a more careful analysis would permit to
discuss also higher values. We confine to codegree $\leq 8$
because $8$ is the smallest value giving rise to the nice example
discussed above.

So, let $n=2$ and set $S=X$. Recall that $c_2(J_1(L))-L^2 =
e(S)+4(g-1)$, where $e(S)$ is the topological Euler--Poincar\'e
characteristic of $S$ and $g=g(L)$. Suppose that $(S,L)$ is
neither $(\Bbb P^2, \Cal O(e))$, $e=1,2$, nor a scroll. Then
$c_2(J_1(L))-L^2 >0$ \cite{LPS1, Prop. A.1}. If $S$ is not
(birationally) ruled, then $e(S) \geq 0$ by the Castelnuovo--De
Franchis theorem. Moreover, $g \geq 2$, with equality if and only
if $S$ is the $K3$ double plane \cite{LP3, Thm. 3.1}, in which
case, however, $e(S)=24$. Thus $e(S)+4(g-1) \geq 8$ if $S$ is
non-ruled. In particular $c_2(J_1(L))
> 8$. Now suppose that $S$ is ruled. Due to our assumptions on
$(S,L)$ we know that $g \geq 1$, and equality occurs if and only
if $S$ is a Del Pezzo surface and $L=-K_S$. For such surfaces we
have $e(S)=12-L^2$ by Noether's formula. Hence $c_2(J_1(L)) =12$.
Assume that $g \geq 2$. If $S \cong \Bbb P^2$, then $g \geq 3$ by
Clebsch formula, hence $e(S)+4(g-1) \geq 3 +8=11$. So,
$c_2(J_1(L))> 12$. On the other hand, if $S \not\cong \Bbb P^2$,
then there exists a birational morphism $\eta:S \to S_0$, where
$S_0$ is a $\Bbb P^1$-bundle over a smooth curve of genus
$q=h^1(\Cal O_S)$. Thus $e(S)=e(S_0) + s = 4(1-q)+s$, where $s$ is
the number of blowing-ups $\eta$ factors though. So we have
$e(S)+4(g-1) = 4(1-q)+s + 4(g-1) \geq 4(g-q).$ As $(S,L)$ is not a
scroll, we know that $K_S+L$ is nef, hence
$$0 \leq (K_S+L)^2 = K_S^2+2(K_S+L)L - L^2 \leq 8(1-q)+4(g-1)-L^2
< 4(1+g-2q).$$ This says that $g \geq 2q$. All cases with $g \leq
1$ being already considered, we conclude that $g-q \geq 1$
equality occurring only for $g=2$. Since $L$ is ample and spanned,
taking into account \cite{LP3, Theorem 3.1} we see that $2=g=q+1$
only for the pair $(S,L)$ in (8.1.3)

So, apart from the pair in $(8.1.3)$ we have $g-q \geq 2$, and
then $e(S)+4(g-1) \geq 8$. In particular, $c_2(J_1(L)) > 8$. The
discussion above proves the following

\proclaim{(8.1.4) Proposition} Let $(S,L, V)$ be as in $(0.0)$,
with $\dim S =2$ and suppose that $(S,V)$ has tame codegree $\leq
8$. If $(S,L)$ is neither $(\Bbb P^2, \Cal O(e))$, $e=1,2$, nor a
scroll, then $\text{\rm{codeg}}(S,V)=8$ and $(S,L,V)$ is as in
$(8.1.3)$.
\endproclaim
Another question is the following:

\noindent{\bf (8.2)} In view of (0.4), for  $(X,L,V)$ as in (0.0)
with $\dim(X) \geq 2$ we can define $R(X,V)=\{D \in |V|:\; D\;
\text{is reducible or non-reduced}\} \subseteq \Cal D(X,V)$. The
following conjecture is quite similar to  \cite{BDL, Conjecture
1}:

\proclaim{(8.2.1) Conjecture} Let $(X,L,V)$ be a triplet as in
$(0.0)$ such that $\dim(X) \geq 2$, $N>n$. Then $R(X,V)=\Cal
D(X,V)$ if and only if $(X,L)$ is either $(\Bbb P^2, \Cal O_{\Bbb
P^2}(2))$ or a scroll over a curve.
\endproclaim

In fact, in \cite{BDL, Conjecture 1} the requirement on the
discriminant locus is replaced with the following condition: for
all $x \in X$,  (i) $|V-2x|\ne \emptyset$ and (ii) any $D \in
|V-2x|$ is reducible or non-reduced. Let us note that (i) is
equivalent to $\dim(|V|) \geq n+1$ and we cannot drop out this
hypothesis, allowing $\phi_V(X)=\Bbb P^N$, as example (d) in (3.3)
shows. A first evidence for this conjecture is that it is true in
the classical case, \cite{BDL, Prop. 7}.

Finally we comment on the following problem.

\noindent{\bf (8.3)} A relevant result in the classical setting is
the so called Landman's parity theorem. The precise statement is
as follows: for $(X,L,V)$ as in (0.0) with $\phi_V$ an embedding
the difference between the defect and the dimension of $X$ is an
even number. This result is not known to be true or false in the
ample and spanned case. First proof of this theorem comes from
Landman \cite{L}, \cite{K, II (22)} and is essentially
topological. In fact since the codimension of the discriminant is
bigger than one one can construct a pencil in $|V|$ not cutting
$\Cal D$, i.e., all its elements are smooth. A consequence of the
existence of these pencils is that the equalities of Betti numbers
of sections of $X$ provided by the Lefschetz theorem go further
for positive defect varieties. This is also true in the ample and
spanned case. Last part of the proof relies on the fact that the
singular locus of a general element of $\Cal D$ is very well known
and provides a vanishing cycle and a monodromy relation giving the
parity result. This last part cannot be applied to the ample and
spanned case. Another proof of Landman's parity theorem can be
found in \cite{Ei}. In the classical setting, when the
discriminant is not a hypersurface, the singular locus of a
general element in $\Cal D$ is a linear space $\Cal T$ of
dimension $k>0$ and any of its points is a non-degenerate
quadratic singularity. Then, the second fundamental form gives a
symmetric isomorphism between the normal bundle $N_{\Cal T/X}$ and
the twist of its dual $N_{\Cal T/X}^\vee(1)$. The symmetry of the
isomorphism and basic considerations on normal bundles have
several relevant consequences like parity theorem (among others).
This puts in relation (8.3) with \cite{LPS1, Conjecture 2.11}, the
conjecture stating that the the singular locus of a general
element in $\Cal D$ is a disjoint union of linear spaces $\Cal T$
of dimension $k$; in particular $X$ is swept out by these linear
spaces (for $\Cal T$ to be linear we mean isomorphic to $\Bbb P^k$
and $\Cal T L^{n-k}=1$). Also \cite{A, Prop. 2.5} shows that
\cite{LPS1, Conjecture 2.11} implies parity theorem.

\vskip0.2cm

\eightpoint {{\it Acknowledgements}. During the preparation of
this paper the first author has been supported by the MIUR of the
Italian Government in the framework of the PRIN "Geometry on
Algebraic Varieties" (Cofin 2004) and by the University of Milano
(FIRST 2003). The second author has been partially supported by
the Project of the Spanish Government BFM2003-03971/MATE.}

\tenpoint

\enddocument